\documentclass{conm-p-l}

\newtheorem{theorem}{Theorem}[section]

\newtheorem{corollary}[theorem]{Corollary}
\newtheorem{conjecture}[theorem]{Conjecture}

\theoremstyle{definition}

\theoremstyle{remark}

\numberwithin{equation}{section}

\input psfig.sty

\def\hf{\frac{1}{2}}

\def\bC{\mathbf C}
\def\bQ{\mathbf Q}

\def\hV{\hat V}
\def\hP{\hat P}
\def\bi{\mathbf i}
\def\boZ{\mathbf Z}
\def\hboZ{\hf\mathbf Z}
\def\bg{\mathbf g}
\def\hg{\hat g}
\def\hW{\widehat W}
\def\ch{\check h}
\def\chel{\check\ell}
\def\bgh{\mathbf{\hat g}}

\def\cC{{\mathcal C}}
\def\cF{{\mathcal F}}

\def\cH{{\mathcal H}}

\def\cO{{\mathcal O}}

\def\cT{{\mathcal T}}

\def\cW{{\mathcal W}}

\def\zd{\zeta_d}

\def\lra{\leftrightarrow}
\def\x{\times}
\def\ox{\otimes}
\def\ni{\noindent}
\def\l{\lambda}
\def\la{\langle}
\def\ra{\rangle}

\begin{document}

\title{Fusion Rules for Affine Kac-Moody Algebras}

\author{Alex J. Feingold}
\address{Dept. of Math. Sci., 
The State University of New York,
Binghamton, New York 13902-6000
}
\email{alex@math.binghamton.edu}

\thanks{I wish to express my thanks to Prof. N. Sthanumoorthy and the
other organizers of the Symposium for inviting me, and for their generous 
hospitality.}

\subjclass[2000]{Primary 17B67, 17B65, 81T40;
Secondary 81R10, 05E10}
\date{}

\keywords{Fusion Rules, Affine Kac-Moody Algebras}

\maketitle

\tableofcontents

\section{Introduction}

Fusion rules play a very important role in conformal field theory \cite{Fu},
in the representation theory of vertex operator algebras \cite{FLM, FHL, FZ}, 
and in quite a few other areas. This paper is not meant to be comprehensive, but 
should be a useful introduction to the subject, with major focus on the 
algorithmic aspects of computing fusion rules in the case of affine Kac-Moody
algebras. I have included many explicit examples and figures illustrating 
the rank 2 cases which can be done graphically on a sheet of paper. The 
Kac-Walton algorithm \cite{Kac, Wal} for fusion coefficients is closely related 
to the Racah-Speiser algorithm for tensor product decompositions, which was the 
subject of my thesis \cite{F1,F2}. I have included here some discussion of this 
relationship and some implications of my thesis for the computation of fusion
coefficients. In Theorems \ref{mynew1} and \ref{mynew2}, for fixed dominant 
integral weights $\lambda$ and $\mu$, I determine the values of level $k$ for 
which all tensor product multiplicities, $Mult_{\lambda,\mu}^\nu$, are equal 
to the corresponding level $k$ fusion coefficients, $N_{\lambda,\mu}^{(k)\ \nu}$, 
for all dominant integral $\nu$. I have recalled the results of Parasarathy, 
Ranga Rao and Varadarajan \cite{PRV} on tensor product multiplicities in 
Theorem \ref{prv} , and the results of Frenkel and Zhu \cite{FZ} on fusion 
coefficients in Theorem \ref{fz}. I have included a conjecture on
fusion coefficients which I believe is a restatement of the Frenkel-Zhu theorem
in a form which shows it to be a beautful generalization of the PRV theorem. 
In joint work \cite{AFW,FW} we have tried to understand fusion
rules from a combinatorial point of view which is quite different from the 
approaches of others \cite{BMW,BKMW,T}. The idea for our new approach was
inspired by our work on explicit spinor constructions \cite{FFR,FRW}.
In \cite{AFW} we explained  
all of the $(p,q)$-minimal model fusion rules \cite{Wa} from elementary
2-groups. The $(p,q)$-minimal models are a certain series of highest weight
representations of the Virasoro algebra \cite{KR} which also have
the structure of a vertex operator algebra \cite{FLM}, and modules
for it \cite{FZ}. In \cite{FW} we explained the fusion rules for all
positive integral levels for type $A_n$ affine Kac-Moody algebras if $n = 1$
or $n = 2$. That work is explained in this paper.   

This paper is an expanded version of two lectures I presented at the 
Ramanujan International Symposium on Kac-Moody Lie Algebras and
Applications, Jan. 28 - 31, 2002, Ramanujan Institute for Advanced
Study in Mathematics, University of Madras, Chennai, India. It was a great
honor to be invited to this symposium, and I was pleased to be able to include
a connection of my work with some work of Ramanujan, whose genius continues to 
inspire great mathematics all around the world. 

\section{Definition of Fusion Algebra}

Let us begin with the definition of fusion algebra given by J. Fuchs \cite{Fu}. 

A fusion algebra $F$ is a finite dimensional commutative associative
algebra over $\bQ$ with some basis 
$$B = \{x_a\ |\ a\in A\}$$ 
so that the structure constants $N_{a,b}^c$ defined by
$$x_a \cdot x_b = \sum_{c\in A} N_{a,b}^c x_c$$
are non-negative integers. There must be a distinguished
index $\Omega\in A$ with the following properties.
Define a matrix
$$C = [C_{a,b}] = [N_{a,b}^\Omega]$$ 
and define an associated ``conjugation'' map $\cC: F\to F$ by
$$\cC(x_a) = \sum_{b\in A} C_{a,b} x_b.$$
It is required that $\cC$ be an involutive automorphism of $F$, so 
$\cC^2 = I_F$ and $C^2 = I$. Because $0\leq N_{a,b}^c \in\boZ$, 
either $C = I$ or $C$ must be an order 2 permutation matrix, that is, there is a
permutation $\sigma:A\to A$ with $\sigma^2 = 1$ and
$$C_{a,b} = \delta_{a,\sigma(b)}.$$
Since $\cC$ is an automorphism, we must also have 
$$\cC(x_a) \cdot \cC(x_b) = \cC(x_a \cdot x_b),$$
that is,
$$x_{\sigma(a)} \cdot x_{\sigma(b)} = \sum_{c\in A} N_{a,b}^c
x_{\sigma(c)}$$
which means that 
$$N_{{\sigma(a)},{\sigma(b)}}^{\sigma(c)} = N_{a,b}^c.$$
Write $\sigma(a) = a^+$ and call $x_{a^+}$ the conjugate of
$x_a$. Use it to define the non-negative integers
$$N_{a,b,c} = N_{a,b}^{c^+}$$
which, by commutativity and associativity of the algebra product, are
completely symmetric in $a$, $b$ and $c$. To see this, note that commutativity 
means $N_{a,b}^{c} = N_{b,a}^{c}$ for all $a,b,c\in A$, so 
$N_{a,b,c} = N_{b,a,c}$. Associativity means
$$\sum_{d\in A} N_{a,b}^{d}\ N_{d,c}^{e} = \sum_{d\in A} N_{b,c}^{d}\ N_{a,d}^{e}$$
for all $a,b,c,e\in A$. Taking $e = \Omega$ and using 
$N_{a,b}^{\Omega} = \delta_{a,b^+} = \delta_{a^+,b}$, this gives 
$N_{a,b}^{c^+} = N_{b,c}^{a^+}$, so $N_{a,b,c} = N_{b,c,a}$. This order 3 
cyclic permutation and the transposition switching $a$ and $b$ generate all 
permutations of $a$, $b$ and $c$. Using this we also find
\begin{equation} \notag
\begin{split}
N_{\Omega,b}^c &= N_{\Omega,b,c^+} = N_{b,c^+,\Omega} = N_{b,c^+}^{\Omega^+} \\
&= N_{b^+,c}^{\Omega} = C_{b^+,c} = \delta_{b,c} \end{split}
\end{equation}
which means $x_\Omega$ is a multiplicative identity element in $F$, so we
write $x_\Omega = 1$. It also follows that $\Omega^+ = \Omega$.

Here are some examples of fusion algebras, which we will later see come
from representations of affine Kac-Moody algebras of some ``level''. 
The algebras are presented by giving a table of products of the basis elements.
These fusion rule tables were produced by the computer program of Bert 
Schellekens, called ``Kac'', available from his webpage:\qquad
http://norma.nikhef.nl/$\sim$t58/

\section{Examples of Fusion Algebras}

\bigskip
\begin{center}
{\bf Table 1:} Fusion Table for $A_1$ of level $k = 2$ \\[2mm]
\framebox{
\begin{tabular}{c|c|c|c}
[i]$\cdot$[j] & [0] & [1] & [2] \\ \hline
[0] & [0] & [1] & [2] \\ \hline
[1] &   & [0] & [2] \\ \hline
[2] &   &   & [0]+[1] 
\end{tabular}
} 
\end{center}
\bigskip

In this example, we have
$$A = \{\Omega = 0,1,2\},\qquad
B = \{x_0 = [0], x_1 = [1], x_2 = [2]\}.$$
From Table 1 we see that, for example,
$$[2]\cdot [2] = 1 [0] + 1 [1] + 0 [2]$$
and we can read off particular structure constants, for example,
$$N_{2,2}^{0} = 1\ \ N_{2,2}^{1} = 1,
\ \ N_{2,2}^{2} = 0.$$
It is also easy to see that
$$C = \bmatrix 1 &0 & 0 \\ 0 &1 &0 \\ 0 &0 &1\endbmatrix.$$
We would get the same fusion table for $B_2$ of level $k=1$.

\bigskip
\begin{center}
{\bf Table 2:} Fusion Table for $A_1$ of level $k = 3$ \\[2mm]
\framebox{
\begin{tabular}{c|c|c|c|c}
[i]$\cdot$[j] & [0] & [1] & [2] & [3]\\ \hline
[0] & [0] & [1] & [2] & [3]\\ \hline
[1] &   & [0] & [3] & [2] \\ \hline
[2] &   &   & [0]+[2] & [1]+[3]  \\ \hline
[3] &   &   &   & [0]+[2]
\end{tabular}
} 
\end{center}
\bigskip

In this example, we have
$$A = \{\Omega = 0,1,2,3\},\qquad 
B = \{[0], [1], [2], [3]\},$$
$$[2]\cdot [3] = 0 [0] + 1 [1] + 0 [2] + 1 [3],$$
$$N_{i,j}^{0} = \delta_{i,j}\ \hbox{so}\ C = I_4.$$

\bigskip
\begin{center}
{\bf Table 3:} Fusion Table for $A_2$ of level $k = 2$ \\[2mm]
\framebox{
\begin{tabular}{c|c|c|c|c|c|c}
[i]$\cdot$[j] & [0] & [1] & [2] & [3] & [4] & [5]\\ \hline
[0] & [0] & [1] & [2] & [3] & [4] & [5]\\ \hline
[1] &   & [2] & [0] & [4] & [5] & [3]\\ \hline
[2] &   &   & [1] & [5] & [3] & [4] \\ \hline
[3] &   &   &   & [0]+[3] & [1]+[4] & [2]+[5] \\ \hline
[4] &   &   &   &   & [2]+[5] & [0]+[3] \\ \hline
[5] &   &   &   &   &   & [1]+[4] 
\end{tabular}
} 
\end{center}
\bigskip

In this example, we have
$$A = \{\Omega = 0,1,2,3,4,5\},\qquad
B = \{[0], [1], [2], [3], [4], [5]\},$$
$$[1]\cdot [2] = [0] = [0]\cdot[0],\qquad [3]\cdot[3] = [0]+[3] = [4]\cdot[5],$$
$$C = \bmatrix 1&0&0&0&0&0 \\ 0&0&1&0&0&0 \\ 0&1&0&0&0&0
\\0&0&0&1&0&0\\ 0&0&0&0&0&1 \\ 0&0&0&0&1&0 \endbmatrix.$$

\bigskip
\begin{center}
{\bf Table 4:} Partial fusion Table for $A_2$ of level $k = 3$ \\[2mm]
\framebox{
\begin{tabular}{c|c|c|c|c}
[i]$\cdot$[j] & [0] & [1] & [2] & [9]\\ \hline
[0] & [0] & [1] & [2] & [9]\\ \hline
[1] &   & [2] & [0] & [9]\\ \hline
[2] &   &   & [1] & [9] \\ \hline
[9] &   &   &   & [0]+[1]+[2]+2[9]
\end{tabular}
} 
\end{center}
\bigskip

In this example, we have
$$A = \{\Omega = 0,1,2,3,4,5,6,7,8,9\},\qquad
B = \{[0], [1], [2],\cdots, [9]\},$$
$$[9]\cdot [9] = [0]+[1]+[2]+2[9].$$
Note that this is the first example where a coefficient exceeds 1: 
$N_{9,9}^{9} = 2$.

Here is the fusion table for the affine algebra of type $B_2$ of level $2$.

\bigskip
\begin{center}
{\bf Table 5:} Fusion Table for $B_2$ of level $k = 2$ \\[2mm]
\framebox{
\begin{tabular}{c|c|c|c|c|c|c}
[i]$\cdot$[j] & [0] & [1] & [2] & [3] & [4] & [5]\\ \hline
[0] & [0] & [1] & [2] & [3] & [4] & [5]\\ \hline
[1] &   & [0] & [3] & [2] & [4] & [5]\\ \hline
[2] &   &   & [0]+[4]+[5] & [1]+[4]+[5] & [2]+[3] & [2]+[3] \\ \hline
[3] &   &   &   & [0]+[4]+[5] & [2]+[3] & [2]+[3] \\ \hline
[4] &   &   &   &   & [0]+[1]+[5] & [4]+[5] \\ \hline
[5] &   &   &   &   &   & [0]+[1]+[4]
\end{tabular}
} 
\end{center}
\bigskip

\section{Notations}

Now we will introduce notations and discuss how fusion algebras are 
associated with representations of untwisted affine Kac-Moody algebras of fixed
level. Let $\bg$ be a finite dimensional simple Lie algebra of rank $N-1$ with
Cartan matrix $A = [a_{ij}]$,
and let 
$$\bgh = \bg\otimes \bC[t,t^{-1}] \oplus \bC c \oplus \bC d$$ 
be the corresponding affine algebra with derivation $d = -t\frac{d}{dt}$ adjoined 
as usual. Let $H$ be the Cartan subalgebra of $\bg$ and let 
$$\cH = H \oplus \bC c \oplus \bC d$$ 
be the Cartan subalgebra of $\bgh$. 
The simple roots and the fundamental weights of $\bg$ are linear functionals
$$\alpha_1,\cdots,\alpha_{N-1} \quad\hbox{and}\quad\lambda_1,\cdots,\lambda_{N-1},$$
respectively, in the dual space $H^*$. Let the integral weight lattice $P$ be 
the $\boZ$-span of the fundamental weights, and let
$$P^+ = \{n_1\lambda_1+\cdots+n_{N-1}\lambda_{N-1}\ |\ 
0\leq n_1,\cdots,n_{N-1}\in\boZ\}$$ 
be the set of dominant integral weights of $\bg$, and let
$$\theta = \sum_{i=1}^{N-1} \ell_i \alpha_i$$
be the highest root of $\bg$. The symmetric bilinear form $(\cdot,\cdot)$ on $H^*$
is determined by 
$$a_{ij} = \langle\alpha_i,\alpha_j\rangle 
= \frac{2(\alpha_i,\alpha_j)}{(\alpha_j,\alpha_j)},\quad 1\leq i,j\leq N-1$$
and the normalization $(\theta,\theta) = 2$. The fundamental weights are 
determined by the conditions $\langle\lambda_i,\alpha_j\rangle = \delta_{ij}$
for $1\leq i,j\leq N-1$, and the special ``Weyl vector''
$$\rho = \sum_{i=1}^{N-1} \lambda_i$$
will play an important role in several formulas. It is useful to define
$${\check\lambda} = \frac{2\lambda}{(\lambda,\lambda)}\quad\hbox{for any }
0\neq \lambda\in H^*,$$
so we can write $(\lambda_i,{\check\alpha}_j) = \delta_{ij}$ and 
$a_{ij} = (\alpha_i,{\check\alpha}_j)$. We may also express
$$\theta = \sum_{i=1}^{N-1} \chel_i {\check\alpha}_i\qquad\hbox{so}\qquad
\chel_i = \frac{\ell_i (\alpha_i,\alpha_i)}{2}.$$
The dual Coxeter number of $\bg$ is defined to be 
$$\ch = 1 + \sum_{i=1}^{N-1} \chel_i  = 1 + \langle\rho,\theta\rangle.$$
The Weyl group $W$ of $\bg$ is defined to be the group of endomorphisms of $H^*$
generated by the simple reflections corresponding to the simple roots, 
$$r_i(\lambda) = \lambda - (\lambda,{\check\alpha}_i),\qquad 1\leq i\leq N-1.$$
This is a finite group of isometries which preserve $P$. 
There is a partial order defined on $H^*$ defined by
$$\lambda \leq \mu \quad\hbox{ if and only if }\quad \mu - \lambda = 
\sum_{i=1}^{N-1} k_i \alpha_i \quad\hbox{for some } 0\leq k_i\in\boZ.$$
For $\lambda\in P^+$ let $V^\lambda$ denote the finite dimensional irreducible
$\bg$-module with highest weight $\lambda$. It has the weight space decomposition
$V^\lambda = \bigoplus_{\beta\in H^*} V^\lambda_\beta$, where 
$$V^\lambda_\beta = \{v\in V^\lambda\ |\ h\cdot v = \beta(h) v, \forall h\in H\}$$
is the $\beta$ weight space of $V^\lambda$. Of course, there are only finitely
many $\beta\in H^*$ such that $V^\lambda_\beta$ is nonzero, and we denote by 
$\Pi^\lambda$ that finite set of such $\beta$. The dual space $(V^\lambda)^* = 
Hom(V^\lambda,\bC)$ is also an irreducible highest weight $\bg$-module, called
the contragredient module of $V^\lambda$. The action of $\bg$ on $(V^\lambda)^*$
is given by 
$$(x\cdot f)(v) = -f(x\cdot v) \quad\hbox{ for }\quad x\in\bg,\ f\in (V^\lambda)^*,
\ v\in V^\lambda.$$
The highest weight of $(V^\lambda)^*$ is denoted by $\lambda^+ = \lambda^*$, and 
equals the negative of the lowest weight of $V^\lambda$. For example, in the case 
when $\bg$ is of type $A_{N-1}$, if $\lambda = \sum_{i=1}^{N-1} n_i \lambda_i$ 
then $\lambda^+ = \sum_{i=1}^{N-1} n_{N-i} \lambda_i$. 

The simple roots and the fundamental weights of 
$\bgh$ are linear functionals  
$$\alpha_0, \alpha_1,\cdots,\alpha_{N-1} \quad\hbox{and}\quad
\Lambda_0,\Lambda_1,\cdots,\Lambda_{N-1},$$
respectively, in the dual space $\cH^*$. The simple roots of $\bg$ form a basis
of $H^*$ (as do the fundamental weights), and we identify them with linear 
functionals in $\cH^*$ having the same values on $H\subseteq\cH$ and being 
zero on $c$ and $d$. Let $c^*$ and $d^*$
in $\cH^*$ be the functionals which are zero on $H$ and which satisfy
$$c^*(c) = 1, \quad c^*(d) = 0, \quad d^*(c) = 0, \quad d^*(d) = 1.$$
Extend the bilinear form $(\cdot,\cdot)$ to $\cH^*$ by letting
$$(c^*,H^*) = 0 = (d^*,H^*),\qquad (c^*,c^*) = 0 = (d^*,d^*),\ \hbox{and}\ 
(c^*,d^*) = 1.$$
Then $\alpha_0 = d^* - \theta$ and 
$$\Lambda_0 = c^*, \quad \Lambda_i 
= \ell_i \frac{(\alpha_i,\alpha_i)}{2}\ c^* + \lambda_i = \chel_i\ c^* + \lambda_i,
\quad 1\leq i\leq N-1,$$
are determined by the conditions $\langle\Lambda_i,\alpha_j\rangle = \delta_{ij}$
for $0\leq i,j\leq N-1$. Let the integral weight lattice $\hP$ be 
the $\boZ$-span of the fundamental weights, and let
$$\hP^+ = \{\sum_{i=0}^{N-1} n_i \Lambda_i \ |\ 0\leq n_i\in\boZ\}$$ 
be the set of dominant integral weights of $\bgh$. 
The affine Weyl group $\hW$ of $\bgh$ is the group of endomorphisms of $\cH^*$
generated by the simple reflections corresponding to the simple roots,
$$r_i(\Lambda) = \Lambda - (\Lambda,{\check\alpha}_i) \alpha_i,\qquad 
0\leq i\leq N-1.$$
This is an infinite group of isometries which preserve $\hP$. 
The canonical central element, $c\in\bgh$ acts on an irreducible $\bgh$-module 
as a scalar $k$, called the level of the module. 
We will only discuss modules with highest weight $\Lambda\in \hP^+$, which are 
the ``nicest'' in that they have affine Weyl group symmetry and satisfy the 
Weyl-Kac character formula. An irreducible highest weight $\bgh$-module is 
uniquely determined by its highest weight 
$$\Lambda = \sum_{i=0}^{N-1} n_i \Lambda_i \in \hP^+$$
and, if we define $\ell_0 = 1 = \chel_0$, then 
$$k = \Lambda(c) = \sum_{i=0}^{N-1} n_i\Lambda_i(c) 
= \sum_{i=0}^{N-1} n_i \ell_i \frac{(\alpha_i,\alpha_i)}{2} =
\sum_{i=0}^{N-1} n_i \chel_i.$$
For fixed $k$ there are only finitely many $\Lambda\in \hP^+$ with 
$\Lambda(c) = k$, and we denote that finite set by $\hP_k^+$. It is easy to 
see that $\hW$ preserves the level $k$ weights 
$\{\Lambda\in\hP\ |\ \Lambda(c) = k\}$. The affine hyperplane determined by 
the condition $\Lambda(c) = k$ can be projected onto $H^*$ and the corresponding 
action of $\hW$ is such that the simple reflections $r_i$ for $1\leq i\leq N-1$ 
act as they were defined originally on $H^*$, as isometries generating 
the finite Weyl group $W$ of $\bg$. But the new affine reflection $r_0$ acts
as $r_0(\lambda) = \lambda - (\lambda,\theta)\theta + k\theta 
= r_{\theta}(\lambda) + k\theta$, the composition of reflection $r_\theta$ 
and the translation by $k\theta$, which is not an isometry on $H^*$. 

Irreducible $\bgh$-modules $\hV^\Lambda$ of level $k \geq 1$ are indexed by 
$\hP_k^+$, but we can also index them by certain weights of $\bg$ as follows. 
From the formulas above we can write
$$\Lambda = \sum_{i=0}^{N-1} n_i \Lambda_i 
= k c^* + \sum_{i=1}^{N-1} n_i \lambda_i.$$
So there is a bijection between $\hP_k^+$ and the set of weights 
$\lambda = \sum_{i=1}^{N-1} n_i \lambda_i$
such that 
$$k = n_0 + \sum_{i=1}^{N-1} n_i \ell_i \frac{(\alpha_i,\alpha_i)}{2} 
= n_0 + \sum_{i=1}^{N-1} n_i \chel_i = n_0 + \langle\lambda,\theta\rangle.$$
Since $n_0 \geq 0$, this is equivalent to the ``level $k$ condition''
$$\langle\lambda,\theta\rangle 
= \sum_{i=1}^{N-1} n_i \chel_i \leq k.$$
Define the set 
$$P_k^+ = \{\lambda = \sum_{i=1}^{N-1} n_i \lambda_i\in P^+\ |\ 
\langle\lambda,\theta\rangle \leq k\}$$
and let the index set $A$ (as in the fusion algebra definition) be $P_k^+$. 
Then we see that irreducible modules on level $k$  
correspond to $N$-tuples of nonnegative integers
$$(n_0,n_1,\cdots,n_{N-1})\ \hbox{such that }
k = \sum_{i=0}^{N-1} n_i \ell_i \frac{(\alpha_i,\alpha_i)}{2} = 
\sum_{i=0}^{N-1} n_i \chel_i .$$ 
Such an $N$-tuple corresponds to
$$\Lambda = k c^* + n_1\lambda_1 + \cdots + n_{N-1} \lambda_{N-1}.$$
Fix level $k\geq 1$ and write the fusion algebra product (which has not been
defined yet!)
$$[\lambda]\cdot [\mu] = \sum_{\nu\in P_k^+} N_{\lambda,\mu}^{\nu}\ [\nu].$$
The distinguished identity element, $[0]$, corresponds to $\Lambda = k c^*$, 
and for each $[\lambda]$ there is a distinguished conjugate $[\lambda^+]$
such that $N_{\lambda,\mu}^{0} = \delta_{\mu,\lambda^+}$. 
Knowing $N_{\lambda,\mu}^{\nu}$
is equivalent to knowing the completely symmetric coefficients
$$N_{\lambda,\mu,\nu} = N_{\lambda,\mu}^{\nu^+}.$$
Let $\cF(\bg,k)$ denote this fusion algebra.

In the case when $\bg = sl_N$ is of type $A_{N-1}$, we have $\ell_i = 1$ 
and $(\alpha_i,\alpha_i) = 2$ for
$0\leq i\leq N-1$, so the set of all weights of level 1,
$$\hP_1^+ = \{\Lambda_i\ |\ 0\leq i\leq N-1\}$$ 
is precisely the set of the fundamental weights of $\bgh$, and
$$P_1^+ = \{0, \lambda_i\ |\ 1\leq i\leq N-1\}.$$
The level 1 fusion algebra $\cF(sl_N,1)$ has a basis 
$\{[0],[1],\cdots,[N-1]\}$ (in Schellekens notation, $[i]$ corresponds to 
$\lambda_i$ for $1\leq i\leq N-1$) and the fusion rules are given by the group 
$\boZ_N$, the weight lattice modulo the root lattice of $\bg$. This means 
$$N_{\lambda,\mu}^\nu = \delta_{\lambda_+\mu,\nu}$$
where the addition takes place in the quotient group of the weight
lattice modulo the root lattice. So $\cF(sl_N,1)$ is the group
algebra $\bQ[\boZ_N]$. 

For $\bg = so(2N)$, $N\geq 4$, of type $D_N$ and rank $N$, 
we have $\ell_i = 1$ for $i = 0,1,N-1,N$, $\ell_i = 2$ for 
$2\leq i\leq N-2$, and $(\alpha_i,\alpha_i) = 2$ for $0\leq i\leq N$,
so the level 1 weights
$$\hP_1^+ = \{\Lambda_i\ |\ i = 0,1,N-1,N\}$$
are the fundamental weights corresponding to the four endpoints of the
affine Dynkin diagram of $\bgh$, and 
$$P_1^+ = \{0, \lambda_i\ |\ i = 1,N-1,N\}$$
where $\lambda_1$ is the highest weight of the natural representation 
(of dimension $2N$), $\lambda_{N-1}$ and $\lambda_{N}$ are highest weights 
of the half-spinor representations (each of dimension $2^{N-1}$. 
The group structure of the weight lattice modulo the root lattice is known
to be the Klein 4-group $\boZ_2 \x \boZ_2$ if $N$ is even, $\boZ_4$ if $N$ 
is odd, and $P_1^+$ is a set of coset representatives for that quotient group
in either case. The fusion algebra for $\bg = so(2N)$ on level 1 is then 
\begin{equation} \notag
\cF(so(2N),1) = \begin{cases}
\bQ[\boZ_2 \x \boZ_2] &\text{if $N$ is even,}\\
\bQ[\boZ_4] &\text{if $N$ is odd.} \end{cases}
\end{equation}

For $\bg = so(2N+1)$, $N\geq 3$, of type $B_N$ and rank $N$, 
we have $\ell_i = 1$ for $i = 0,1$, $\ell_i = 2$ for $2\leq i\leq N$,
$(\alpha_i,\alpha_i) = 2$ for $0\leq i\leq N-1$ and 
$(\alpha_N,\alpha_N) = 1$. In the special case of $B_2$, we have 
$(\ell_0, \ell_1, \ell_2) = (1,2,1)$, and 
$(\alpha_i,\alpha_i) = 2,1,2$ for $i = 0,1,2$, respectively. 
So for $N\geq 2$, the level 1 weights are
$$\hP_1^+ = \{\Lambda_i\ |\ i = 0,1,N\}.$$
The weight lattice modulo the root lattice is $\boZ_2$. 

For $\bg = sp(2N)$, $N\geq 2$, of type $C_N$ and rank $N$, 
we have $\ell_i = 1$ and $(\alpha_i,\alpha_i) = 2$ for $i = 0,N$,
and $\ell_i = 2$ and $(\alpha_i,\alpha_i) = 1$ for $1\leq i\leq N-1$.  
The weight lattice modulo the root lattice is $\boZ_2$. 
The special case of $N=2$ gives $B_2$. 

We will not give further details about the exceptional algebras, but for
$\bg$ of type $G_2$ we should mention that 
$(\ell_0, \ell_1, \ell_2) = (1,3,2)$, where $\alpha_2$ is the long root
and $(\alpha_1,\alpha_1) = 2/3$.

\section{Algorithms For Tensor Product Decompositions}

There is a close relationship between the product in fusion algebras
associated with an affine Kac-Moody algebra $\bgh$ and tensor product 
decompositions of irreducible $\bg$-modules. 
Let $V^\lambda$ be the irreducible finite dimensional $\bg$-submodule
of $\hV^\Lambda$ generated by a highest weight vector. In the special
case when $\Lambda = k\Lambda_0 = kc$, that finite dimensional 
$\bg$-module is $V^0$, the one dimensional trivial $\bg$-module. 
Since $\bg$ is semisimple, any finite dimensional $\bg$-module is 
completely reducible. Therefore, we can write the tensor product 
of irreducible $\bg$-modules 
$$V^\lambda\otimes V^\mu = \sum_{\nu\in P^+} 
Mult_{\lambda,\mu}^\nu V^\nu$$
as the direct sum of irreducible $\bg$-modules, including multiplicities. 
This decomposition is independent of the level $k$ and is part of the 
basic representation theory of $\bg$. The fusion products
$[\lambda]\cdot [\mu]$ are obtained by a subtle truncation 
of the above summation. 

The Racah-Speiser algorithm gives the formula 
$$Mult_{\lambda,\mu}^\nu = \sum_{w\in W} \epsilon(w) 
Mult_\lambda(w(\nu+\rho) - \mu - \rho)$$
where $W$ is the Weyl group of $\bg$, 
$\epsilon(w) = (-1)^{length(w)}$ is the sign of $w$, 
the Weyl vector $\rho = \sum \lambda_i$ is the sum of 
the fundamental weights of $\bg$, and 
$Mult_\lambda(\beta) = dim(V^\lambda_\beta)$ is the inner multiplicity of 
the weight $\beta$ in $V^\lambda$. Recall that 
$\Pi^\lambda = \{\beta\in H^*\ |\ dim(V^\lambda_\beta) > 0\}$ denotes the set 
of all weights of $V^\lambda$. 

In fact, the only weights $\nu$ for which $Mult_{\lambda,\mu}^\nu$ 
may be nonzero are those of the form $\nu = \beta + \mu$ where 
$\beta\in \Pi^\lambda$. This means the formula is a geometrical algorithm: 

(1) Shift the weight diagram of $V^\lambda$ by adding $\mu + \rho$. 

(2) Use the Weyl group to move all shifted weights into
the dominant chamber, where they accumulate as an alternating sum of
inner multiplicities of $V^\lambda$, adding if the required $w$ is even,
subtracting if it is odd. 

(3) The resulting pattern of numbers will be non-negative
integers, zero if the shifted weight is on a chamber wall, and after shifting
the pattern back by subtracting $\rho$, you will have the ``outer'' tensor 
product multiplicities. 

This algorithm assumes that you can already produce the weight diagram of
any irreducible module, $V^\lambda$, so we should have discussed that first,
but in fact the special case of the Racah-Speiser algorithm when $\mu = 0$
gives a recursion for the inner multiplicities of $V^\lambda$. Since $V^0$
is the trivial one-dimensional module, $V^\lambda\otimes V^0 = V^\lambda$, 
so $Mult_{\lambda,0}^\nu = \delta_{\lambda,\nu}$ and therefore
$$0 = \sum_{w\in W} \epsilon(w) Mult_\lambda(w(\nu+\rho) - \rho)$$
for $\nu\neq\lambda$. One knows that $Mult_\lambda(w\lambda) = 1$  and 
$Mult_\lambda(w\nu) = Mult_\lambda(\nu)$ for all
$w\in W$, so the above formula implies that 
$$Mult_\lambda(\nu) = - \sum_{1\neq w\in W} \epsilon(w) 
Mult_\lambda(\nu+\rho - w\rho)$$
for $\nu\neq\lambda$. Since $\rho > w\rho$ in the partial ordering on weights,
this gives an effective recursion for $Mult_\lambda(\nu)$. It is instructive
to carry out these recursions by hand in the rank 2 cases, where the geometry
is simple to see on a sheet of paper. I have included in the appendices pages 
of type $A_2$, $B_2$ and $G_2$ weight lattices, including the reflecting axes,
and pages with just the reflecting axes. If you make a copy of the former, 
you can put on it the weight diagram of a single irreducible module, $V^\lambda$,
by starting with one dominant weight (make a heavy dot) at position $\lambda$. 
Then find all the dominant weights less than $\lambda$ in the partial ordering.
Apply the Weyl group to that set of weights to get all the weights of the module.
An example of this for type $A_2$ with $\lambda = 3\lambda_1+2\lambda_2$ is given
in the Appendix, Figure 8. To use the Racah recursion formula, make a copy of 
the reflection axes only on a transparency. (Choose the appropriate axis for
the type of algebra from Figures 4 - 6.) Then place the 
weight diagram you made under the transparency, shifted by $\rho$. 
Using the shifted reflecting lines you can see the points
which will be involved in the alternating sum for a given dominant $\nu$ in
the diagram, and find the multiplicity of the $\nu$ weight space. Mark those
multiplicities next to each weight, using the Weyl group action on the 
unshifted weight diagram to mark nondominant weights. 
Now you can use that marked weight diagram to compute the tensor product
of that module with any other by the Racah-Speiser algorithm. You only need to
put the diagram under the transparency of reflecting axes shifted by $\mu+\rho$
and follow steps (2) and (3) above. For example, using the weight diagram in
Figure 8 to compute the tensor product decomposition of $V^\lambda\otimes V^\mu$
for $\lambda = 3\lambda_1+2\lambda_2$ and $\mu = \lambda_1$, one would see the
shifted weight diagram shown in Figure 9, and find that the Racah-Speiser
algorithm gives the answer 
$$V^{4\lambda_1+2\lambda_2} \oplus V^{3\lambda_1+\lambda_2} \oplus 
V^{2\lambda_1+3\lambda_2}.$$
But in that case, it would have been wiser to shift the weight diagram of 
$V^{\lambda_1}$ by $\mu = 3\lambda_1+2\lambda_2$ plus $\rho$ as shown in Figure
16. (Ignore for now the affine reflection line shown there.) The three
weights of that fundamental module each have multiplicity 1, and after shifting
by $4\lambda_1+3\lambda_2$, all of them are strictly inside the dominant 
chamber, so there are no cancellations and each of them gives a highest weight
module in the tensor product decomposition as shown above. 

Another method of recursively computing the weight multiplicities 
$Mult_\lambda(\nu)$ is as follows. Place a clear transparency over the weight 
lattice, locate the weights $w\rho = w(\lambda_1 + \lambda_2)$ for each $w\in W$, 
and make an open circle around each such point, large enough to see the underlying
weight in the diagram. Since the differences between those points and the fixed
point $\rho$ is $\rho - w\rho$, if you rotate the transparency 180 degrees and
place the point $\rho$ over any weight $\nu$ of a weight diagram for $V^\lambda$,
the other open circles of the transparency will lie over the points $\nu + \rho -
w\rho$, which will be strictly above $\nu$ in the partial ordering of weights. 
The Racah recursion formula can then be implemented by taking the alternating sum 
of the multiplicities of those circled weights, assumed to have been already found 
by the initial data $Mult_\lambda(w\lambda) = 1$, or by the application of Weyl 
group symmetry $Mult_\lambda(w\beta) = Mult_\lambda(\beta)$ to multiplicities 
already found recursively. I have combined in Figure 7 the diagrams of the Weyl 
conjugates of $\rho$ for each type. If you copy this page onto a transparency,
it can be used as described above to recursively compute weight multiplicities
of irreducible modules for any of the rank 2 algebras. After the page is rotated 
by 180 degrees, the open circle corresponding to $\rho$ should be placed over
the weight to be computed. It will be the alternating sum of the weights under
the other circles, where the plus or minus signs inscribed in the circles indicate
whether to add or subtract. It is well known that for type $A_2$ the resulting
pattern of multiplicities is easy to describe. The weight diagrams for type 
$A_2$ consist of concentric hexagonal shells, which may degenerate into triangles
towards the center. The outer shell consists of weights all of whose multiplicities
are equal to $1$. The weights on the next hexagonal shell inward have multiplicity
$2$, and each successive shell inward has all multiplicities one more than the
one outside it. This pattern continues until the hexagonal shell becomes a triangle.
The multiplicity of each weight on that triangle, and on all weights further inward,
is the same, one more than the multiplicity on that innermost hexagon. For example,
in Figure 8, the weight diagram consists of two hexagonal shells and one triangular
shell. The outer hexagonal shell has 15 weights, each with multiplicity equal 
to $1$, the next hexagonal shell has 9 weights each with multiplicity equal to $2$,
and the inner triangular shell has 3 weights each with multiplicity $3$. As a check 
on this, note that $42 = (15)(1) + (9)(2) + (3)(3)$ is then the dimension of the 
irreducible $A_2$-module in Figure 8. Using $n_1 = 3$ and $n_2 = 2$, this agrees 
with the formula 
$$\dim(V^\lambda) = (n_1+n_2+2)(n_1+1)(n_2+1)/2$$
for an irreducible $A_2$-module $V^\lambda$ with 
$\lambda = n_1\lambda_1+n_2\lambda_2$.

In my thesis \cite{F1,F2} I studied certain patterns which occur in the tensor 
product decomposition of a fixed irreducible $\bg$-module, $V^\lambda$, with all 
other modules $V^\mu$. For fixed $\lambda$, as $\mu$ varies there are only a finite 
number of different patterns of outer multiplicities which can occur, and there
are sets of values for $\mu$ for which the pattern is constant. I called those
zones of stability for tensor product decompositions, and they can be understood
from the geometrical point of view of the Racah-Speiser algorithm.
If the weight diagram of $V^\lambda$ is shifted parallel to one of the fundamental
weights, say by $\mu + m\lambda_i$, there is a least value $m_i$ such that for 
$m\geq m_i$, the set of shifted weights, $\Pi^\lambda + \mu + m\lambda_i + \rho$ 
is contained in the union of the images of the fundamental chamber under $W(i)$, 
the subgroup of the Weyl group generated by the simple reflections $r_j$, $j\neq i$. 
These are the chambers containing the weights $k\lambda_i$ for $k\geq 1$. 
If $m$ exceeds $m_i$, the only $w\in W$ which may make nonzero contributions to
the outer multiplicity are those from $W(i)$, and those fix $\lambda_i$. The
geometrical reflection process which generates the tensor product multiplicities
is therefore the same for each $m\geq m_i$.
While the highest weights of the modules occurring increase by the number of
$\lambda_i$'s added, their outer multiplicities stay constant. In fact, we have
the following precise result from \cite{F2} about when a particular weight 
$\beta$ of $V^\lambda$, reaches the zone of stability.

\medskip
\begin{theorem}\label{mythesis} Let $\lambda,\mu\in P^+$ and $\beta\in\Pi^\lambda$ be
such that $\beta + \mu\in P^+$. Let 
$$\beta - r_{\beta,j}\alpha_j,\cdots,\beta,\cdots\beta + q_{\beta,j}\alpha_j$$
be the $\alpha_j$ weight string through $\beta$. If
$\langle\mu,\alpha_j\rangle \geq q_{\beta,j}$ then 
$$Mult_{\lambda,\mu}^{\beta+\mu} = 
Mult_{\lambda,\mu+\lambda_j}^{\beta+\mu+\lambda_j}.$$ \end{theorem}

Since $\langle\mu+\lambda_j,\alpha_j\rangle = \langle\mu,\alpha_j\rangle + 1$,
it is clear that $\langle\mu,\alpha_j\rangle \geq q_{\beta,j}$ implies
$$Mult_{\lambda,\mu}^{\beta+\mu} = 
Mult_{\lambda,\mu+m\lambda_j}^{\beta+\mu+m\lambda_j}\qquad\hbox{for all }m\geq 1.$$
This result shows that for fixed $\lambda\in P^+$ and fixed
$\beta\in\Pi^\lambda$, the tensor product multiplicities 
$Mult_{\lambda,\mu}^{\beta+\mu}$ have zones of stability as $\mu$ varies, and 
it is sufficient to study the finite number of $\mu$ such that 
$\langle\mu,\alpha_j\rangle \leq q_{\beta,j}$ for $1\leq j\leq N-1$.

For example, using the weight diagram of $V^\lambda = V^{3\lambda_1+2\lambda_2}$ 
for $A_2$ shown in Figure 8, look at the weight $\beta = 2\lambda_1+\lambda_2$. 
The $\alpha_1$ weight string through this $\beta$ goes from 
$-4\lambda_1+4\lambda_2 = \beta - 3\alpha_1$ to $4\lambda_1 = \beta + \alpha_1$, 
so $q_{\beta,1} = 1$. The $\alpha_2$ weight string through this $\beta$ goes from
$4\lambda_1 - 3\lambda_2 = \beta - 2\alpha_2$ to $\lambda_1 + 3\lambda_2 = 
\beta + \alpha_2$ so $q_{\beta,2} = 1$. Theorem \ref{mythesis} then says that if 
$\mu = n_1 \lambda_1 + n_2 \lambda_2$ then
$\langle\mu,\alpha_1\rangle = n_1 \geq q_{\beta,1} = 1$ implies 
$$Mult_{\lambda,\mu}^{\beta+\mu} = 
Mult_{\lambda,\mu+\lambda_1}^{\beta+\mu+\lambda_1}$$
and $\langle\mu,\alpha_2\rangle = n_2 \geq q_{\beta,2} = 1$ implies
$$Mult_{\lambda,\mu}^{\beta+\mu} = 
Mult_{\lambda,\mu+\lambda_2}^{\beta+\mu+\lambda_2}.$$
In Figure 9 we can see the weight $\beta$, with multiplicity 2, shifted by 
$\mu+\rho = 2\lambda_1 + \lambda_2$, in position for the reflection 
process, which will reduce it by 1 because of the weight $r_2(\beta+\mu+\rho)$.
Since $\mu = m\lambda_1$ for $m\geq 1$ satisfies the conditions of Theorem  
\ref{mythesis} for $\alpha_1$, we have 
$$1 = Mult_{\lambda,m\lambda_1}^{\beta+m\lambda_1}\qquad\hbox{for all } m\geq 1.$$
It is clear that as $m$ increases, the reflection process
yields the same result as $\beta+\mu+\rho$ shifts further along the line parallel 
to $\lambda_1$. In contrast, $\mu = \lambda_1$ does not satisfy the conditions
of Theorem \ref{mythesis} for $\alpha_2$ and we can see that adding $\lambda_2$ 
to $\mu$ means shifting the weight diagram in Figure 9 by $\lambda_2$, which 
leads to a different reflection process for the shifted $\beta$ and a different 
multiplicity. 

There is another important result about tensor product coefficients which 
played a role in my thesis. I will always be grateful to Prof. Bertram Kostant
for drawing my attention to the following beautiful result of Parthasarathy, 
Ranga Rao and Varadarajan \cite{PRV}, which I have rewritten in the form I
found most useful in my thesis.

\medskip
\begin{theorem}\label{prv} \cite{PRV} Let $\lambda,\mu\in P^+$ and 
$\beta\in\Pi^\lambda$ be such that $\beta + \mu\in P^+$. Let $\ell = rank(\bg)$ 
and let $e_j\in\bg$ be a root vector corresponding
to the simple root $\alpha_j$ for $1\leq j\leq \ell$. Then 
$$Mult_{\lambda,\mu}^{\beta+\mu} = 
dim\{v\in V^\lambda_\beta\ |\ e_j^{\langle\mu,\alpha_j\rangle+1} v = 0, 
1\leq j\leq \ell\}.$$\end{theorem}

\section{Algorithms For Fusion Product Coefficients}

Let $N_{\lambda,\mu}^{(k)\ \nu}$ denote the fusion product coefficient
at level $k$. Then the Kac-Walton algorithm (\cite{Kac}, p. 288, \cite{Wal}) 
expresses this as an alternating
sum of tensor product multiplicities:
$$N_{\lambda,\mu}^{(k)\ \nu} = \sum_{w\in \hW} \epsilon(w) 
Mult_{\lambda,\mu}^{w(\nu+\rho)-\rho}$$
where $\hW$ is the affine Weyl group acting on the weight lattice of $\bg$
with the action of the simple reflections of $W$ as usual, but with 
$$r_0(\beta) = r_\theta(\beta) + (k+{\ch})\theta.$$
Here $r_\theta$ is reflection with respect to the highest root $\theta$ of $\bg$,
and $\ch$ is the dual Coxeter number of $\bg$. In the case when
$\bg = sl_N$, $\ch = N$, $W$ is the symmetric group $S_N$, and
$\theta = \sum \alpha_i$ is the sum of the simple roots of $\bg$. Let
$T_y(x) = x+y$ be the function which translates by vector $y$. Then it is 
easy to see that 
$$T_{s\theta}\ r_\theta\ T_{-s\theta}(\beta) = r_\theta(\beta) + 2s\theta$$
which will equal $r_0(\beta)$ if $s = \frac{1}{2}(k+{\ch})$. Therefore, 
$r_0$ is reflection with respect to the shifted hyperplane
perpendicular to $\theta$, translated by $\frac{1}{2}(k+{\ch})\theta$. 

Using the Racah-Speiser formula in the Kac-Walton formula gives
a formula for fusion coefficients as an alternating sum of inner
multiplicities:
$$N_{\lambda,\mu}^{(k)\ \nu} = \sum_{w\in \hW} \epsilon(w) 
Mult_\lambda(w(\nu+\rho) - \mu - \rho)$$
which has a nice geometrical interpretation as before, but using the
affine Weyl group $\hW$ instead of $W$. 

The only weights $\nu$ for which $N_{\lambda,\mu}^{(k)\ \nu}$ 
may be nonzero are those of the form 
$\nu = \beta + \mu$ where $\beta\in \Pi^\lambda$. 
The geometrical interpretation of this formula is now as follows: 

(1) Shift the weight diagram of $V^\lambda$ by adding $\mu + \rho$. 

(2) Use the affine Weyl group to move all shifted weights into the part of
the dominant chamber bounded by the reflection wall of the affine reflection,
$r_0$, where they accumulate as an alternating sum of
inner multiplicities of $V^\lambda$, adding if the required $w\in \hW$ 
is even, subtracting if it is odd. 

(3) The resulting pattern of numbers will be non-negative
integers, zero if the shifted weight is on a reflection wall, and after 
shifting the pattern back by subtracting $\rho$, you will have the fusion 
product coefficients. 

To get a better intuitive understanding of this algorithm, it is useful to do
some rank 2 cases using the diagrams from the Appendix. To include the new
affine reflection, $r_0$, you need to make a transparency for the reflection 
line corresponding to the highest root, $\theta$, and you need to know the 
dual Coxeter number, $\ch$, where
\begin{align} \notag
\theta &=  
\begin{cases} \lambda_1 + \lambda_2, &\text{in the $A_2$ case;} \\
2\lambda_1 &\text{in the $B_2$ case;} \\ 
\lambda_2  &\text{in the $G_2$ case.} 
\end{cases} &
\ch &=  
\begin{cases} 3, &\text{in the $A_2$ case;} \\
3 &\text{in the $B_2$ case;} \\ 
4  &\text{in the $G_2$ case.} 
\end{cases} 
\end{align}
For $\bg$ of type $A_2$, level 2, the following table gives the correspondence 
between the fusion algebra labels [i] used in Table 3, the triples 
$(n_0,n_1,n_2)$ whose sum equals the level, and the weights $\lambda = 
n_1 \lambda_1 + n_2 \lambda_2\in P_2^+$:

\bigskip
\begin{center}
{\bf Table 6:} Label-Weight Correspondence for $A_2$ of level $k = 2$ \\[2mm]
\framebox{
\begin{tabular}{c|c|c}
[i] & $(n_0,n_1,n_2)$ & $\lambda = n_1 \lambda_1 + n_2 \lambda_2$\\ \hline
[0] & $(2,0,0)$ &  $0$  \\ \hline
[1] & $(0,2,0)$ &  $2\lambda_1$  \\ \hline
[2] & $(0,0,2)$ &  $2\lambda_2$  \\ \hline
[3] & $(0,1,1)$ &  $\lambda_1+\lambda_2$   \\ \hline
[4] & $(1,0,1)$ &  $\lambda_2$   \\ \hline
[5] & $(1,1,0)$ &  $\lambda_1$  
\end{tabular}
} 
\end{center}
\bigskip

To check, for example, the fusion product $[3]\cdot[3] = [0] + [3]$ from
Table 3, we would take the weight diagram of $V^{\lambda_1+\lambda_2}$,
the adjoint representation, and shift it by $\mu + \rho = 
2\lambda_1+2\lambda_2$, and use the affine Weyl group to move all shifted
weights into the part of the dominant chamber bounded by the affine reflecting
line. See Figure 10 and verify that after shifting back by $\rho$ the surviving
highest weights are $0$ and $\lambda_1+\lambda_2$, each with multiplicity 1. 
(Note that the tensor product multiplicity of $\lambda_1+\lambda_2$ would
have been 2, but the affine reflection line reduced it by one, and killed two
other weights which were on it.)

In the next two tables, for $\bg$ of type $B_2$, levels 1 and 2, respectively, 
we give the correspondence between the fusion algebra labels [i] used in 
Tables 1 and 5, respectively, the triples 
$(n_0,n_1,n_2)$ whose sum equals the level, and the weights 
$\lambda = n_1 \lambda_1 + n_2 \lambda_2$:

\bigskip
\begin{center}
{\bf Table 7:} Label-Weight Correspondence for $B_2$ of level $k = 1$ \\[2mm]
\framebox{
\begin{tabular}{c|c|c}
[i] & $(n_0,n_1,n_2)$ & $\lambda = n_1 \lambda_1 + n_2 \lambda_2$\\ \hline
[0] & $(1,0,0)$ &  $0$  \\ \hline
[1] & $(0,0,1)$ &  $\lambda_2$ \\ \hline
[2] & $(0,1,0)$ &  $\lambda_1$ 
\end{tabular}
} 
\end{center}

\bigskip
\begin{center}
{\bf Table 8:} Label-Weight Correspondence for $B_2$ of level $k = 2$ \\[2mm]
\framebox{
\begin{tabular}{c|c|c}
[i] & $(n_0,n_1,n_2)$ & $\lambda = n_1 \lambda_1 + n_2 \lambda_2$\\ \hline
[0] & $(2,0,0)$ &  $0$  \\ \hline
[1] & $(0,0,2)$ &  $2\lambda_2$  \\ \hline
[2] & $(1,1,0)$ &  $\lambda_1$  \\ \hline
[3] & $(0,1,1)$ &  $\lambda_1+\lambda_2$   \\ \hline
[4] & $(0,2,0)$ &  $2\lambda_1$   \\ \hline
[5] & $(1,0,1)$ &  $\lambda_2$  
\end{tabular}
} 
\end{center}
\bigskip

We may check the level 1 fusion products from Table 1 as follows. For 
$[1]\cdot[1] = [0]$, take the weight diagram of $V^{\lambda_2}$, shift it by 
$\lambda_2 + \rho = \lambda_1+2\lambda_2$, and use the affine Weyl group to 
move all shifted weights into the part of the dominant chamber bounded by the 
affine reflecting line. See Figure 11 and verify that after shifting back by 
$\rho$ the only surviving highest weight is $0$. For $[1]\cdot[2] = [2]$, shift 
the diagram of $V^{\lambda_2}$ by $\lambda_1 + \rho = 2\lambda_1+\lambda_2$, and 
after the same process (see Figure 12) find the only surviving highest 
weight is $\lambda_1$. For $[2]\cdot[2] = [0]+[1]$, shift the weight
diagram of $V^{\lambda_1}$ by $\lambda_1 + \rho = 2\lambda_1+\lambda_2$, and
after reflecting (see Figure 13) find the only surviving highest weights are 
$0$ and $\lambda_2$. 

In Figure 14 check the level 2 fusion product $[3]\cdot[3] = [0]+[4]+[5]$ from 
Table 5 by shifting the weight diagram of $V^{\lambda_1+\lambda_2}$  by 
$\lambda_1+\lambda_2 + \rho = 2\lambda_1+2\lambda_2$, and after the affine
reflection process (with the affine reflection line located as it should be for
level 2) verify that (after shifting back by $\rho$) the surviving highest weights
are $0$, $2\lambda_1$ and $\lambda_2$.

In comparing the Kac-Walton algorithm with the one for $Mult_{\lambda,\mu}^\nu$, 
we see that the shifting is the same, and all reflections coming from $w\in W$ are
the same, but there are more contributions from the extra elements in $\hW$. 
Elements of $W$ are sufficient to reflect all weights of the diagram into
the dominant chamber, but some may be on the side of the reflection wall of
$r_0$ not containing the origin. One application of $r_0$ would then move
the weight to the other side, but perhaps take it out of the dominant chamber,
requiring more reflections from $W$ to move it back into the dominant chamber.
For example, in Figure 15, which is just Figure 9 with the affine
reflection line added in the position appropriate for level 5, we see that 
the shifted highest weight falls on the affine reflection line. So the tensor
product decomposition given in the last section is truncated by removing that
highest weight on the line to give the corresponding fusion product. This can 
be seen more clearly in Figure 16, where the smaller module $V^{\lambda_1}$ is
shifted by $3\lambda_1+2\lambda_2+\rho$. In Figure 17 we see an example where
many weights of the shifted module are on the far side of the affine reflection 
line, and where application of $r_0$ does not bring the weight into $P_k^+$. 

We would like to briefly discuss how the result in Theorem \ref{mythesis} 
on zones of stability for tensor product multiplicities might 
give some information about such zones for fusion coefficients. 
It is clear that $P_k^+ \subseteq P_{k+1}^+$, and for fixed $\lambda$ and $\mu$, 
increasing level $k$ means that the affine reflection wall will move further 
away from the origin. For $\alpha$ any root, the reflection with respect to the 
hyperplane perpendicular to $\alpha$ is
$$r_\alpha(\lambda) = \lambda - \langle\lambda,\alpha\rangle\ \alpha$$
and if we write $T_\beta(\lambda) = \lambda + \beta$ for translation by $\beta$,
then
$$r_0 r_\theta(\lambda) = \lambda + (k+\ch)\theta = T_{(k+\ch)\theta}(\lambda)$$
is a translation in $\hW$. It is easy to check the relation
$r_\alpha T_\beta r_\alpha = T_{r_\alpha(\beta)}$, which implies that
$$w T_{(k+\ch)\theta} w^{-1} = T_{(k+\ch)w(\theta)}$$
is a translation in $\hW$ for each $w\in W$. Since $\theta$ is the highest root
of $\bg$, it's orbit under $W$ is the set of all long roots. 
So for each long root, $\beta$, the translation $T_{(k+\ch)\beta}$ and its 
inverse are even elements of $\hW$. These translations generate an abelian 
subgroup $\cT$ of $\hW$ and the relation
$w T_\beta w^{-1} = T_{w(\beta)}$ shows that $\hW$ is the 
semi-direct product of $W$ and $\cT$. A fundamental domain for the action of 
$\hW$ on the set of all weights $P$ can be determined by writing any element of
$\hW$ as a translation from $\cT$ followed by an element of $W$. The translations
$T_{(k+\ch)n\beta}$, $n\in\boZ$, allow any weight $\lambda$ to be moved to a weight 
$\mu$ such that $-(k+\ch) \leq \langle\mu,\beta\rangle \leq (k+\ch)$. These 
inequalities say that $\mu$ is between the shifted hyperplanes fixed by
$T_{(k+\ch)\beta} r_\beta$ and by $r_\beta T_{(k+\ch)\beta}$. Doing this for
each positive long root $\beta$ allows us to move $\lambda$ to a weight $\mu$ 
in the closure of a fundamental domain for $\cT$, 
$$\cF_k = \{\mu\in P\ |\ -(k+\ch) \leq \langle\mu,w\theta\rangle \leq (k+\ch), 
\forall w\in W\},$$
the region bounded by all such pairs of shifted hyperplanes. That region is 
obviously $W$-invariant, and each weight in it can be moved by $W$ into the
dominant chamber, $P^+$. So a fundamental domain for $\hW$ would be the
intersection $\cF_k \cap P^+ = P_{k+\ch}^+$ and  
$$\cF_k = W(P_{k+\ch}^+) = \bigcup_{w\in W} w(P_{k+\ch}^+).$$
Let
$$\cF_k' = \{\mu\in P\ |\ -(k+\ch) < \langle\mu,w\theta\rangle < (k+\ch), 
\forall w\in W\}$$
be the interior of $\cF_k$. Then $\cF_k'$ is also $W$-invariant, and for any
translation $T\in\cT$, if $\cF_k'\cap T(\cF_k')$ is nonempty then $T$ is the
identity element. The boundary walls of $\cF_k$ expand as k increases so there 
is a minimum value of $k$ for which $\Pi^\lambda + \mu + \rho \subseteq \cF_k'$.
For example, in Figure 9 we see that the shifted weight diagram is contained in
the interior of the large hexagon, $\cF_6'$.

\medskip
\begin{theorem}\label{mynew1} For $\lambda,\mu\in P^+$, if $k$ is large enough so
that $\Pi^\lambda +\mu+\rho\subseteq\cF_k'$, then for all $\nu\in P_k^+$ we have
$$N_{\lambda,\mu}^{(k)\ \nu} = Mult_{\lambda,\mu}^{\nu}.$$ \end{theorem}

\noindent{\bf Proof:} For $\lambda,\mu\in P^+$ we see that if 
$\Pi^\lambda +\mu+\rho\subseteq\cF_k'$ then only elements of $W$ can bring 
those shifted weights into $P^+_{k+\ch}$, and none
go on the fixed hyperplane of $r_0$. When using the Kac-Walton algorithm to 
compute the fusion coefficients in the product $[\lambda]\cdot[\mu]$ this 
condition guarantees that the only nonzero contributions may come from affine 
Weyl group elements which are actually in $W$, matching the expression in the 
Racah-Speiser algorithm for the tensor product coefficients and giving
the equality of the fusion and tensor product coefficients as claimed.
$\square$\medskip

I have the following result for finding that minimum value of $k$ which 
makes the above happen.

\medskip
\begin{theorem}\label{mynew2} For $\lambda,\mu\in P^+$, we have
$$\langle\lambda+\mu,\theta\rangle \leq k \hbox{ if and only if }
\ \Pi^\lambda +\mu+\rho\subseteq\cF_k'.$$ \end{theorem}

\noindent{\bf Proof:} First note that $\langle\lambda+\mu,\theta\rangle \leq k$
is equivalent to $\langle\lambda+\mu+\rho,\theta\rangle < k + \ch$ because
$\langle\lambda+\mu+\rho,\theta\rangle = (\lambda+\mu+\rho,\theta) = 
(\lambda+\mu,\theta) + (\rho,\theta) = \langle\lambda+\mu,\theta\rangle+\ch-1$.
Since $\lambda\in \Pi^\lambda$, if $\Pi^\lambda +\mu+\rho\subseteq\cF_k'$ 
then $\lambda+\mu+\rho \in \cF_k'$, so 
$\langle\lambda+\mu+\rho,\theta\rangle < k + \ch$, so we get
$\langle\lambda+\mu,\theta\rangle \leq k$. 

Now suppose that we have the above inequality. To show the containment we break 
the argument into two steps. We will show 
$$(1)\quad \Pi^{\lambda+\mu+\rho} \subseteq \cF_k'\qquad 
\hbox{and then we will show}\qquad
(2)\quad \Pi^\lambda +\mu+\rho \subseteq \Pi^{\lambda+\mu+\rho}.$$
Both $\Pi^{\lambda+\mu+\rho}$ and $\cF_k'$ are $W$-invariant sets, so each
consists of the $W$-conjugates of their dominant integral elements. So if 
$\Pi^{\lambda+\mu+\rho}\cap P^+ \subseteq \cF_k'\cap P^+$ then we get (1).
For any $\beta\in\Pi^{\lambda+\mu+\rho}\cap P^+$, we know 
$\beta \leq \lambda+\mu+\rho$, so $\beta = \lambda+\mu+\rho -
\sum_{i=1}^{N-1} k_i \alpha_i$ with $0\leq k_i\in\boZ$. Then we have
$$\langle\beta,\theta\rangle = \langle\lambda+\mu+\rho,\theta\rangle - 
\sum_{i=1}^{N-1} k_i \langle\alpha_i,\theta\rangle.$$
We also know that $\langle\alpha_i,\theta\rangle = (\alpha_i,\theta) =
(\theta,\alpha_i) = \langle\theta,\alpha_i\rangle\ (\alpha_i,\alpha_i)/2$, 
but $(\alpha_i,\alpha_i) > 0$ and $\theta\in P^+$ since it is the highest
weight of the adjoint representation, so $\langle\theta,\alpha_i\rangle \geq 0$.
Therefore, 
$$\langle\beta,\theta\rangle\leq\langle\lambda+\mu+\rho,\theta\rangle < k+\ch$$
so $\beta\in \cF_k'\cap P^+$. 

It is well-known that for any $\lambda\in P^+$, 
$$\Pi^\lambda = \{w\beta\ |\ \beta\in P^+, \beta\leq\lambda, w\in W\}.$$
Let $\beta\in\Pi^\lambda$ so $\beta\leq \lambda$ and 
$\beta+\mu+\rho\leq \lambda+\mu+\rho$. For some $w\in W$ we have 
$w(\beta+\mu+\rho)\in P^+$. We know $\mu \geq w\mu$ and $\rho\geq w\rho$
since $\mu,\rho\in P^+$, and $\lambda\geq w\beta$ since $w\beta\in\Pi^\lambda$.
Then 
$$\lambda+\mu+\rho - w(\beta+\mu+\rho) = (\lambda - w\beta) + 
(\mu - w\mu) + (\rho - w\rho) \geq 0.$$
But $\lambda+\mu+\rho \geq w(\beta+\mu+\rho) \in P^+$ means 
$w(\beta+\mu+\rho) \in \Pi^{\lambda+\mu+\rho}$ so 
$\beta+\mu+\rho \in \Pi^{\lambda+\mu+\rho}$ by the $W$-invariance of
$\Pi^{\lambda+\mu+\rho}$. Note that this proof of (2) does not use the 
inequality involving $k$. $\square$
\medskip

In \cite{FZ} the following formula for fusion coefficients for affine algebras 
was obtained using the theory of vertex operator algebras. (Also see \cite{GW}.)

\medskip
\begin{theorem}\label{fz} \cite{FZ} Let $\lambda,\mu,\nu\in P^+_k$, and
let $e_\theta$ be a root vector of $\bg$ in the $\theta$ root space of
$\bg$. Let $v_\l^\l \in V^\l$ be a highest weight vector. Then
the level $k$ fusion coefficient $N_{\l,\mu,\nu}^{(k)}$ equals 
the dimension of the vector space $T_k(\l,\mu,\nu) = $
$$\{f\in Hom_\bg(V^\l \ox V^\mu \ox V^\nu,\bC)\ |\ 
f(v_\l^\l \ox e_\theta^{k-\la\l,\theta\ra+1} v^\mu \ox
v^\nu) = 0,\ \forall v^\mu\in V^\mu,\ \forall v^\nu\in V^\nu \}.$$
\end{theorem}
\medskip

It is clear that the $k$-dependent condition on $f$ in $T_k(\l,\mu,\nu)$
will be trivially satisfied for any $\nu$ when the operator 
$e_\theta^{k-\la\l,\theta\ra+1}$ is the zero operator on $V^\mu$,
and in that case $dim(T_k(\l,\mu,\nu)) = Mult^0_{\l,\mu,\nu} = 
Mult^{\nu^+}_{\l,\mu}$ equals the multiplicity of the trivial module in
the triple tensor product $V^\l \ox V^\mu \ox V^\nu$ which equals the
multiplicity of the contragrediant module $V^{\nu^+}$ in $V^\l \ox V^\mu$.
Consider the decomposition of $V^\mu$ into irreducible $sl_2$-modules with
respect to the subalgebra $\bg_\theta \subseteq \bg$ with basis $e_\theta$, 
$f_\theta$ in the $-\theta$ root space of $\bg$ and $h_\theta = 
[e_\theta,f_\theta]$. It is well-known that any finite dimensional irreducible 
representation $V(n)$ of $sl_2$ is uniquely determined by it's highest eigenvalue 
for $h_\theta$, $0\leq n \in \boZ$, and that $dim(V(n)) = n+1$. Using the
well-known action of $e_\theta$ on $V(n)$, it is easy to see that 
$e_\theta^{n+1}$ is the zero operator on $V(n)$. In the decomposition of $V^\mu$
into $sl_2$-modules, there is a component $V(n)$ with largest $n$, and so
$e_\theta^{n+1}$ is the zero operator on that and all other components. 
It is not hard to see that the largest $n$ is $\la\mu,\theta\ra$, which 
corresponds to the $\bg_\theta$-submodule generated by the highest weight vector 
of $V^\mu$. Then the combined results of Theorems \ref{mynew1} and \ref{mynew2} 
follow from Theorem \ref{fz} because the condition on $k$ which guarantees 
equality of fusion and tensor coefficients is that 
$$k-\la\l,\theta\ra + 1 \geq \la\mu,\theta\ra + 1.$$
It is interesting to see how the geometrical aspects of the Kac-Walton and 
Racah-Speiser algorithms give this same result. 

If we do not demand equality of fusion and tensor product coefficients for all 
weights of the shifted weight diagram, we can still get a condition which 
guarantees it for a fixed weight of $\Pi^\l$. 

For each $\beta\in\Pi^\l$ such that $\beta+\mu\in P^+$, 
there is a minimum value of $k$, denoted by
$k_{max} = k_{max}(\beta,\lambda,\mu)$, such that for any $w\in\hW$, 
$w(\beta+\mu+\rho)\in\Pi^\lambda + \mu+\rho$ implies $w\in W$. Assuming that
$k\geq \langle\lambda,\theta\rangle$ and $k\geq \langle\mu,\theta\rangle$ so
that $\lambda,\mu\in P^+_k$, and that $k\geq \langle\beta+\mu,\theta\rangle$
so that $\beta+\mu\in P^+_k$, 
if $k\geq k_{max}$ then the discussion above shows that 
$$N_{\lambda,\mu}^{(k)\ (\beta+\mu)} = Mult_{\lambda,\mu}^{\beta+\mu}.$$

\medskip
\begin{conjecture} \label{AJFConj1} For $\lambda,\mu\in P^+$, $\beta\in\Pi^\lambda$
and $k$ large enough so that $\lambda,\mu,\beta+\mu\in P^+_k$, suppose that
$r_0(\beta+\mu+\rho) \notin \Pi^\l + \mu + \rho$. Then for any $w\in\hW$,
we have
$$w(\beta+\mu+\rho)\in\Pi^\l + \mu+\rho\qquad 
\hbox{implies} \qquad w\in W.$$ \end{conjecture}
\medskip

I re-discovered the following conjecture, which appeared in \cite{Wal2} without
proof. (Thanks to Mark Walton for informing me of his paper after seeing an
earlier version of this paper on the internet arXiv.) As far as I know, it 
remains unproven, but will be the subject of a subsequent publication if I
can prove it. 

\begin{conjecture} \label{AJFConj2} For $\l,\mu\in P^+_k$, $\beta\in\Pi^\l$
such that $\beta+\mu\in P^+_k$, we have $N_{\l,\mu}^{(k)\ (\beta+\mu)}$
equals the dimension of the space 
$$F_k^+(\l,\beta,\mu) = 
\{v\in V^\l_\beta\ |\ e_j^{\la\mu,\alpha_j\ra+1} v = 0, 
1\leq j\leq\ell, \hbox{ and } e_\theta^{k-\la\beta+\mu,\theta\ra+1} v = 0\}.$$ 
\end{conjecture} 

This conjecture is a blending of the PRV and FZ theorems, showing that the
FZ theorem is actually a beautiful generalization of the PRV theorem. It 
implies the following result, which tells the level $k$ at which the fusion
coefficient associated with a single weight $\beta\in\Pi^\l$ equals 
the tensor product multiplicity associated with that weight.

\begin{corollary} \label{AJFCor} Suppose $\l,\mu\in P^+_k$, and
$\beta\in\Pi^\lambda$ is such that $\beta+\mu\in P^+_k$. Let the $\theta$ 
weight string through $\beta$ in $\Pi^\lambda$ be 
$\beta - r\theta,\cdots,\beta,\cdots,\beta + q\theta$. Then 
$k\geq \la\mu,\theta\ra + r$ implies 
$N_{\l,\mu}^{(k)\ (\beta+\mu)} = Mult_{\l,\mu}^{\beta+\mu}$.
\end{corollary}

\noindent{\bf Proof:} If the condition 
$e_\theta^{k-\la\beta+\mu,\theta\ra+1} v = 0$ is satisfied for all
$v\in V^\l_\beta$ then $F_k^+(\l,\beta,\mu) = V^+(\l,\beta,\mu)$
whose dimension is the tensor product multiplicity $Mult_{\l,\mu}^{\beta+\mu}$.
But that condition will be satisfied when $k-\la\beta+\mu,\theta\ra+1 > q$
because that many applications of the operator $e_\theta$ will move $v$
just beyond the $\beta + q\theta$ weight space in the string. We know that 
$r = r_{\l,\beta}$ and $q = q_{\l,\beta}$ depend on $\l$ and on
$\beta$, and satisfy $r - q = \la\beta,\theta\ra$, so the inequality above is
equivalent to $k\geq \la\mu,\theta\ra + r$. 
\medskip

In this approach to fusion coefficients, for fixed values of $\lambda$, $\mu$ and
$\nu$, as the level $k$ varies, we try to determine for what level $k_{max}$ they
reach their maximum, the tensor product coefficient. This should be compared to 
the use of ``threshold levels'' in \cite{BKMW}. The spaces $T_k(\lambda,\mu,\nu)$
for fixed weights as $k$ increases form a filtration of the largest such space,
which is when $k \geq k_{max}(\beta,\lambda,\mu)$, $\nu^+ = \beta+\mu$. There is 
also a $k_{min} = k_{min}(\beta,\lambda,\mu)$ such that 
$N_{\lambda,\mu}^{(k)\ (\beta+\mu)} = 0$ for $k < k_{min}$ but
$N_{\lambda,\mu}^{(k)\ (\beta+\mu)} \neq 0$ for $k = k_{min}$. We may choose 
a basis $B_k$ of each space $T_k$, $k_{min} \leq k \leq k_{max}$ so that each 
$B_k$ is an extension of $B_{k-1}$. Then for each basis vector, $v\in B_{k_{max}}$ 
there is a smallest $k_v$ such that the vector is in $B_{k_v}$, and that $k_v$ 
is called the threshold level of $v$. Knowing the list of all threshold levels
is equivalent to knowing all the fusion coefficients as $k$ varies, but since
there is no canonical choice of basis in the spaces $T_k$, it seems more 
natural to focus on the dimensions of the spaces $T_k$. 

Finally, the tensor product multiplicity
$Mult_{\lambda,\mu}^\nu$ may be part of a zone of uniform decomposition, and
equal to another one with $\mu$ and $\nu$ reduced. For example, in Figure 18
we see the weight diagram from Figure 15 shifted by an additional $\lambda_1$,
which changes the $W$-reflection process for many of the weights, but not for
the weight $2\lambda_1+\lambda_2$ which is in its zone of stability along the
$\lambda_1$ line. But when doing the Kac-Walton algorithm, the affine reflection
cancels that tensor product coefficient because of its $r_0$ symmetry with the
shifted highest weight. But if the level $k$ is increased, then that symmetry
is broken and that fusion coefficient remains constant for all $k\geq 6$.

\section{A Different Approach}

Let $G = \boZ_N^k$. The symmetric group $S_k$ acts on $G$ by 
permuting the $k$-tuples. For $a\in G$, let $[a]$ be the 
orbit of $a$ and $\cO$ be the set of all orbits. 
These orbits are precisely the subsets
$$P(i_0,i_1,\hdots,i_{N-1}) = 
\{x\in\boZ_N^k\ |\ j\hbox{ occurs exactly }
i_j\hbox{ times in }x, 0\leq j\leq N-1\}$$
where $(i_0,i_1,\cdots, i_{N-1})$ is any $N$-tuple of nonnegative integers 
such that 
$$i_0+i_1+\cdots + i_{N-1} = k.$$ 
We now have a bijection between $\cO$ and the set $P_k^+$ when $\bg$ is of type
$A_{N-1}$. For $[a],[b],[c]\in \cO$ we believe the fusion coefficients
$N_{[a],[b]}^{[c]}$ have a combinatorial description in terms of the group $G$. 
The conjugate of $[c]$ is $[-c]$ and we prefer to study the totally
symmetric coefficients 
$$N_{[a],[b],[c]} = N_{[a],[b]}^{[-c]}.$$

We consider the following combinatorial question.  
For $[a],[b],[c]\in \cO$, the group $S_k$ acts on
$$T([a],[b],[c]) = 
\{(x,y,z)\in [a]\x [b]\x [c] \ |\ x + y + z = 0 \}$$
which decomposes into a finite number of orbits under that action. 
Let the number of such orbits be denoted by $M([a],[b],[c])$. 
Determine $M([a],[b],[c])$ and show how it is related to $N_{[a],[b],[c]}$.
For $N = 2,3$ we have the following results.
\medskip

\begin{theorem}\label{fusion1} For $N = 2$, for any integral level $k\geq 1$, with 
notation as above, we have 
$$M([a],[b],[c]) = N_{[a],[b],[c]}.$$ \end{theorem}

\begin{theorem}\label{fusion2} For $N = 3$, for any integral level $k\geq 1$, with 
notation as above, we have 
$$M([a],[b],[c]) = \binom{N_{[a],[b],[c]} + 1}{2}.$$ \end{theorem}
\medskip

In previous work with F. Akman \cite{AFW}, we introduced the 
idea of covering a fusion algebra by a 
finite abelian group and proved that the $(p,q)$-minimal model fusion
algebra, which comes from the discrete series of $0 < c < 1$ 
representations of the Virasoro algebra, can be covered by the group
$Z_2^{p+q-5}$. The basic idea, which is only set up to handle fusion
algebras whose fusion coefficients $N_{ij}^k$ are in $\{0,1\}$, is 
as follows.

\medskip

\noindent{\bf Definition.} Let $(G,+,0)$ be a finite abelian group and let
$G = P_0 \cup P_1 \cup ... \cup P_{N-1}$ be a partition into
$N$ disjoint subsets with $P_0 = \{0\}$.
Let $W$ be an $N$-dimensional vector space over $\bQ$ with basis
$P = \{P_0,P_1,...,P_{N-1}\}$ and define a bilinear 
multiplication on $W$ by the formula
$$P_i * P_j = \sum_{k\in T(i,j)} P_k$$
where 
$$T(i,j) = \{k\ |\ \exists a\in P_i,\exists b\in P_j,\ a+b\in P_k\}.$$
We say that such a partition is associative if the product $*$ is
associative. We say that a group $G$ covers a fusion algebra $F$ if there is 
an associative partition $P$ of $G$ and a bijection
$\Phi$ between $A$ and $P$ which gives an algebra isomorphism between $F$
and $W$ such that $\Phi(\Omega) = P_0$. 

\medskip

As an example of a nontrivial fusion algebra which can be covered by a
group, let the $\cW$ algebra coming from the coset construction of
$$\frac{SU(N)_r \ox SU(N)_s}{SU(N)_{r+s}}$$
be denoted by $\cW_N(r,s)$. The fusion rules for $\cW_3(1,1)$ are as follows.

\begin{center}
{\bf Table 9:} Fusion rules for $\cW_3(1,1)$ \\[2mm]
\renewcommand{\baselinestretch}{1.2}\tiny\normalsize
\framebox{
\begin{tabular}{c|c|c|c|c|c|c}
 [a]$\x$[b] & [0] & [1] & [2] & [3] & [4] & [5] \\[1mm] \hline
 [0] & [0] & [1] & [2] & [3] & [4] & [5] \\[1mm] \hline
 [1] &  & [0]+[1] & [3] & [2]+[3] & [5] & [4]+[5] \\[1mm] \hline
 [2] &   &   & [4] & [5] & [0] & [1] \\[1mm] \hline
 [3] &   &   &   & [4]+[5] & [1] & [0]+[1] \\[1mm] \hline
 [4] &   &   &   &   & [2] & [3] \\[1mm] \hline
 [5] &   &   &   &   &   & [2]+[3]
\end{tabular}
} 
\renewcommand{\baselinestretch}{1}\tiny\normalsize
\end{center}
\bigskip
   
Note that $\{[0], [2], [4]\}$ forms a subgroup isomorphic to $\boZ_3$.
We find that $\boZ_3^2$ covers these fusion rules as follows:
$$
\begin{array}{rrr}
\{(0,0)\} \lra  [0] \quad &  \{(1,2),(2,1)\}  \lra   [1] \quad & \{(1,1)\}  \lra  [2] \\[2mm]
\{(0,2),(2,0)\} \lra  [3] \quad &  \{(2,2)\}  \lra  [4] \quad &  \{(1,0),(0,1)\}  \lra  [5] 
\end{array}
$$

\medskip
   
Let $\bg = sl_2$.
There are $k+1$ $\hg$-modules $\hV^a$ of level $k \geq 1$, 
indexed by spin $a\in \hboZ$ with $0\leq a \leq \frac{k}{2}$. 
We have the tensor product decomposition
$$V^a\otimes V^b = \sum_{|a-b|\leq c\leq a+b} V^c$$
where the sum is only taken over those $c\in \hboZ$ such that
$a+b+c\in\boZ$.
The fusion rules for level $k$ are a simple truncation of that summation: 
\begin{equation} \notag
N_{a,b}^c = 
\begin{cases} 1, &\text{if $|a-b|\leq c\leq a+b$, $a+b+c\in\boZ$, $a+b+c\leq k$;}\\
0, &\text{otherwise.} \end{cases}
\end{equation}
The conditions above imply that $c\leq \frac{k}{2}$.

Alternative way: Re-index the modules $\hV^a$ on level $k$ by 
$m = 2a+1\in \boZ$ with $1\leq m \leq k+1$. Then $m = dim(V^a)$ and we 
write $\hV^a = \hV(2a+1) = \hV(m)$. Let $p = k+2$.
\medskip

\noindent{\bf Definition.} For integer $p\geq 2$ the triple of integers
$(m,m',m'')$ is {\bf p-admissible} when  $0<m,m',m''<p$, 
the sum $m+m'+m''<2p$ is odd, and the ``triangle'' inequalities 
$$m<m'+m'',\quad m'<m+m'', \quad m''<m+m'$$
are satisfied.

\medskip

Then the level $k$ $sl_2$ fusion rules are: 
$N_{m,m'}^{m''} = 1$ if $(m,m',m'')$ is p-admissible, $N_{m,m'}^{m''} = 0$ 
otherwise.

\medskip

\begin{theorem}\label{fusion3} The level $k$ fusion rules for $\bg = sl_2$ define
a fusion algebra $F$ with 
$$A = \{m\in \boZ\ |\ 1\leq m \leq k+1\},$$
distinguished element $\Omega = 1$ and the conjugate $m^+ = m$.
$F$ is covered by the elementary abelian 2-group 
$G = \boZ_2^{k}$ with partition given by
$$P_i = \{g\in G\ |\ \hbox{exactly }i\hbox{ coordinates of }g\hbox{ are
}1\}$$ 
\ni for $0\leq i\leq k$. \end{theorem}
\medskip

The following tables illustrate how the fusion tables for $A_1$
on levels 2 and 3 are covered. 

\begin{center}
{\bf Table 10:} Group $\boZ_2^2$ covering the Fusion Table 
for $A_1$ of level $k = 2$ \\[2mm] 
\renewcommand{\baselinestretch}{1.3}\tiny\normalsize
\framebox{
\begin{tabular}{c|c|c|c}
[i]$\cdot$[j] & (0,0) & (1,1) & (1,0),(0,1)
\\[1mm] \hline
(0,0) & (0,0) & (1,1) & (1,0),(0,1)
\\[1mm] \hline
(1,1) &  & (0,0) & (0,1),(1,0)
\\[1mm] \hline
(1,0),(0,1) &  &  & $\begin{array}{c} (0,0),(1,1) \\ (1,1),(0,0) \end{array}$
\end{tabular}
} 
\renewcommand{\baselinestretch}{1}\tiny\normalsize
\end{center}
\bigskip

\begin{center}
{\bf Table 11:} Group $\boZ_2^3$ covering the Fusion Table 
for $A_1$ of level $k = 3$ \\[2mm]
\setlength{\tabcolsep}{2pt}
\renewcommand{\baselinestretch}{1.3}\tiny\normalsize
\framebox{
\begin{tabular}{c|c|c|c|c}
[i]$\x$[j] & (0,0,0) & (1,1,1) & (0,1,1),(1,0,1),(1,1,0) & (1,0,0),(0,1,0),(0,0,1)
\\[1mm] \hline
 (0,0,0) & (0,0,0) & (1,1,1) & (0,1,1),(1,0,1),(1,1,0) & (1,0,0),(0,1,0),(0,0,1)
\\[1mm] \hline
 (1,1,1) &  & (0,0,0) & (1,0,0),(0,1,0),(0,0,1) & (0,1,1),(1,0,1),(1,1,0)
\\[1mm] \hline
$\begin{array}{c} (0,1,1) \\ (1,0,1) \\ (1,1,0)  \end{array}$
 & & & $\begin{array}{c}  (0,0,0),(1,1,0),(1,0,1) \\ (1,1,0),(0,0,0),(0,1,1)
\\ (1,0,1),(0,1,1),(0,0,0) \end{array}$
& $\begin{array}{c}  (1,1,1),(0,0,1),(0,1,0) \\ (0,0,1),(1,1,1),(1,0,0)
\\ (1,0,0),(1,0,0),(1,1,1) \end{array}$
\\[1mm] \hline
$\begin{array}{c} (1,0,0) \\ (0,1,0) \\ (0,0,1) \end{array}$
&&&&
$\begin{array}{c} 
(0,0,0),(1,1,0),(1,0,1) \\ (1,1,0),(0,0,0),(0,1,1) \\
(1,0,1),(0,1,1),(0,0,0) \end{array}$
\end{tabular}
} 
\renewcommand{\baselinestretch}{1}\tiny\normalsize
\end{center}
\bigskip

The following table illustrates how the fusion table for $A_2$ of level 2
can be covered. In this case all fusion coefficients are $0$ or $1$, but
in order to make this idea work for higher levels of $A_2$, where the
fusion coefficients can be greater than $1$, the method of Theorem 
\ref{fusion2} must be used. 

\begin{center}
{\bf Table 12:} Group $\boZ_3^2$ covering the Fusion Table 
for $A_2$ of level $k = 2$ \\[2mm]
\setlength{\tabcolsep}{2pt}
\renewcommand{\baselinestretch}{1.3}\tiny\normalsize
\framebox{
\begin{tabular}{c|c|c|c|c|c|c}
 [a]$\x$[b]
 & (0,0) & (1,1) & (2,2)
 & $\begin{array}{c} (1,2) \\ (2,1) \end{array}$
 & $\begin{array}{c} (2,0) \\ (0,2) \end{array}$
 & $\begin{array}{c} (1,0) \\ (0,1) \end{array}$
\\[1mm] \hline
 (0,0)
& (0,0) & (1,1) & (2,2)
 & $\begin{array}{c} (1,2) \\ (2,1) \end{array}$
 & $\begin{array}{c} (2,0) \\ (0,2) \end{array}$
 & $\begin{array}{c} (1,0) \\ (0,1) \end{array}$
\\[1mm] \hline
 (1,1)
 &  & (2,2) & (0,0)
 & $\begin{array}{c} (2,0) \\ (0,2) \end{array}$
 & $\begin{array}{c} (0,1) \\ (1,0) \end{array}$
 & $\begin{array}{c} (2,1) \\ (1,2) \end{array}$
\\[1mm] \hline
 (2,2)
 &  &  & (1,1)
 & $\begin{array}{c} (0,1) \\ (1,0) \end{array}$
 & $\begin{array}{c} (1,2) \\ (2,1) \end{array}$
 & $\begin{array}{c} (0,2) \\ (2,0) \end{array}$
\\[1mm] \hline
 $\begin{array}{c} (1,2) \\ (2,1) \end{array}$
 &  &  & 
 & $\begin{array}{c} (2,1),(0,0) \\ (0,0),(1,2) \end{array}$
 & $\begin{array}{c} (0,2),(1,1) \\ (1,1),(2,0) \end{array}$
 & $\begin{array}{c} (2,2),(0,1) \\ (1,0),(2,2) \end{array}$
\\[1mm] \hline
 $\begin{array}{c} (2,0) \\ (0,2) \end{array}$
 &  &  &  & 
 & $\begin{array}{c} (1,0),(2,2) \\ (2,2),(0,1) \end{array}$
 & $\begin{array}{c} (0,0),(1,2) \\ (2,1),(0,0) \end{array}$
\\[1mm] \hline
 $\begin{array}{c} (1,0) \\ (0,1) \end{array}$
 &  &  &  &  & 
 & $\begin{array}{c} (2,0),(1,1) \\ (1,1),(0,2) \end{array}$
\end{tabular}
} 
\renewcommand{\baselinestretch}{1}\tiny\normalsize
\end{center}
\bigskip

I would like to conclude this section with a discussion of the fusion table for 
$A_2$ of level 3, and explain how Table 4 comes from Theorem \ref{fusion2}. 
We must look at the orbits of $\boZ_3^3$ under the symmetric 
group $S_3$. In the notation used by Schellekens, the primaries $[i]$ for 
$0\leq i\leq 9$ correspond to triples $(i_2,i_1,i_0)$ with sum $i_0+i_1+i_2 = 3$, 
as follows:

$$\begin{array}{ccccc}
\!\!\![0] \lra (0,0,3) \ \ & [1] \lra (0,3,0) \ \ & [2] \lra (3,0,0) \ \ & 
[3] \lra (0,1,2) \ \ & [4] \lra (1,2,0) \\[2mm]
\!\!\![5] \lra (2,0,1) \ \ & [6] \lra (0,2,1) \ \ & [7] \lra (2,1,0) \ \ &
[8] \lra (1,0,2) \ \ & [9] \lra (1,1,1) \end{array}$$

Each triple $(i_2,i_1,i_0)$ corresponds to an orbit in $\boZ_3^3$ consisting
of those triples with $i_2$ 2's, $i_1$ 1's and $i_0$ 0's. In particular, for
the primaries given in the earlier partial table, we have

\begin{align} \notag
[0] &\lra (0,0,3) \lra \hbox{ orbit of } (0,0,0)\in\boZ_3^3 \\
[1] &\lra (0,3,0) \lra \hbox{ orbit of } (1,1,1)\in\boZ_3^3 \notag \\
[2] &\lra (3,0,0) \lra \hbox{ orbit of } (2,2,2)\in\boZ_3^3 \notag \\
[9] &\lra (1,1,1) \lra \hbox{ orbit of } (0,1,2)\in\boZ_3^3 \notag 
\end{align}

Denoting the $S_3$-orbit of $(i_2,i_1,i_0)\in\boZ_3^3$ by $[(i_2,i_1,i_0)]$,
we see that each of the orbits
$$[(0,0,0)] = \{(0,0,0)\}, \quad [(1,1,1)] = \{(1,1,1)\}, \quad [(2,2,2)] = \{(2,2,2)\}$$
consists of only one element of $\boZ_3^3$, but the orbit
$$[(0,1,2)] = \{(0,1,2), (0,2,1), (1,0,2), (1,2,0), (2,0,1), (2,1,0)\}$$
consists of the six distinct elements. Looking at the set $T([a],[b],[c])$ when
$[a]$, $[b]$ and $[c]$ are singleton orbits chosen from among $[0]$, $[1]$ 
and $[2]$, is the same as looking at just one equation, 
$$(i,i,i) + (j,j,j) + (k,k,k) = (0,0,0)$$
which has one solution when $i+j+k=0$ in $\boZ_3^3$, none otherwise. This
corresponds to the $\boZ_3$ subtable generated just by $[0]$, $[1]$ and $[2]$.
To understand the rest of the table, it might be easier to look at the $S_3$
orbits of equations $x+y=z$ for $x\in [a]$, $x\in [b]$ and $z\in [c]$. 
When $[a] = [(i,i,i)]$ is a singleton and $[b] = [9] = [(0,1,2)]$, we are 
looking at equations of the form
$$(i,i,i) + (r,s,t) = (r',s',t')$$
where $r,s,t$ are all distinct, which gives $r'=r+i$, $s'=s+i$ and $t'=t+i$ 
are also all distinct, so $[(r',s',t')] = [9]$. There is only one $S_3$-orbit 
of such equations for fixed $[a]$ and $[b]$ so Theorem \ref{fusion2} says that 
$1 = M([a],[9],[9]) = N_{[a],[9],[9]}(N_{[a],[9],[9]}+1)/2$
so $N_{[a],[9],[9]} = 1$, and $0 = M([a],[9],[c])$ for $[c]\neq[9]$. 
Finally, to look at $[9]\cdot[9]$ we must look at equations of the form
$$(r,s,t) + (r',s',t') = (i,j,k)$$
where $r,s,t$ are all distinct and $r',s',t'$ are all distinct. The 36 
possibilities for $(i,j,k)$ include $(0,0,0)$ six times, $(1,1,1)$ six
times, $(2,2,2)$ six times, and each of the six elements of $[9]$ occurs
3 times. One finds that each group of six with $i=j=k$ is a single $S_3$-orbit,
but the 18 equations with $i$,$j$,$k$ distinct fall into three orbits, 
the diagonal orbit of $(0,1,2) + (0,1,2) = (0,2,1)$, the orbit of
$(0,1,2) + (1,2,0) = (1,0,2)$, and the orbit of $(1,2,0) + (0,1,2) = (1,0,2)$.
Theorem \ref{fusion2} then says 
$1 = M([9],[9],[a]) = N_{[9],[9],[a]}(N_{[9],[9],[a]}+1)/2$
so $N_{[9],[9],[a]} = 1$ for $[a] = [0], [1], [2]$, and
$3 = M([9],[9],[9]) = N_{[9],[9],[9]}(N_{[9],[9],[9]}+1)/2$
so $N_{[9],[9],[9]} = 2$. 

\section{A Connection With Ramanujan}

While studying the orbits of $S_k$ acting on $\boZ_n^k$ (in collaboration with 
Michael Weiner and Matthias Beck) we noticed the following. Let 
$$A(n,k,r) = \{(a_1,\cdots,a_k)\in\boZ_n^k\ |\ a_1+\cdots+a_k\equiv r \mod\ n\}$$
and let $M(n,k,r)$ be the number of orbits of $A(n,k,r)$ under the action 
of $S_k$. 
Equivalently, we can represent each such orbit uniquely by a $k$--tuple of 
integers $(a_1,\cdots,a_k)$ where
$0\leq a_j\leq n-1$ for $0\leq j \leq k$, and 
$$a_1\geq a_2\geq\cdots\geq a_k.$$ 
Each such $k$--tuple corresponds to a partition of $a_1+\cdots+a_k$ into
at most $k$ parts, each of which is at most $n-1$. 
Hence, if we denote by $p(a,b,t)$ the number of 
partitions of $t$ into at most $b$ parts, each of which is at most $a$, 
we get the alternative description 
$$M(n,k,r) = \sum_{t \geq 0} p(n-1,k,r+nt).$$
(Here we understand that $p(a,b,0) = 1$ and that $0\leq r < n$.) 
We started to study sums of this type, and proved the beautiful formula 
$$M(n,k,r) = \frac{1}{n+k}\sum_{d|g} \binom{\frac{n+k}{d}}{\frac{n}{d}}\ 
c_d(r),$$
where $g = \gcd(n,k)$ and the sum is over the positive divisors of $g$. 
Here $c_d(r)$ denotes the {\bf Ramanujan sum}, defined for integers $d$ 
and $r$, $d>0$, as 
$$c_d(r) = \sum_{\substack{m=0\\ \gcd(m,d)=1}}^{d-1} \zd^{mr},$$
where $\zd = e^{2\pi\bi/d}$.
One immediately gets the symmetry
$$M(n,k,r) = M(k,n,r).$$
But after posting our results on the internet archives, we learned that
this result was already published in 1999 by Elashvili, Jibladze and Pataraia
\cite{EJP}. Further research in the literature led me back to a 1902 paper by 
von Sterneck \cite{vonS}, who studied partitions into distinct parts, and to
Bachmann \cite{Bac} (Vol. 2, 222--241), who also obtained in 1910 a recursive 
formula for the number of partitions with repetitions allowed, and then to 
Ramanathan \cite{Ram}, who found in 1944 the role of Ramanujan sums in these 
formulas, but did not obtain the beautiful symmetry above. I mentioned this 
history because it shows the far reaching influence of the great mathematician, 
Ramanujan, whose ideas continue to affect the development of mathematics, and 
in whose name we meet at this International Symposium.

\newpage

\section{Appendices}
\centerline{Figure 1: Weight Lattice of Type $A_2$}
\vskip 20pt
\centerline{\psfig{figure=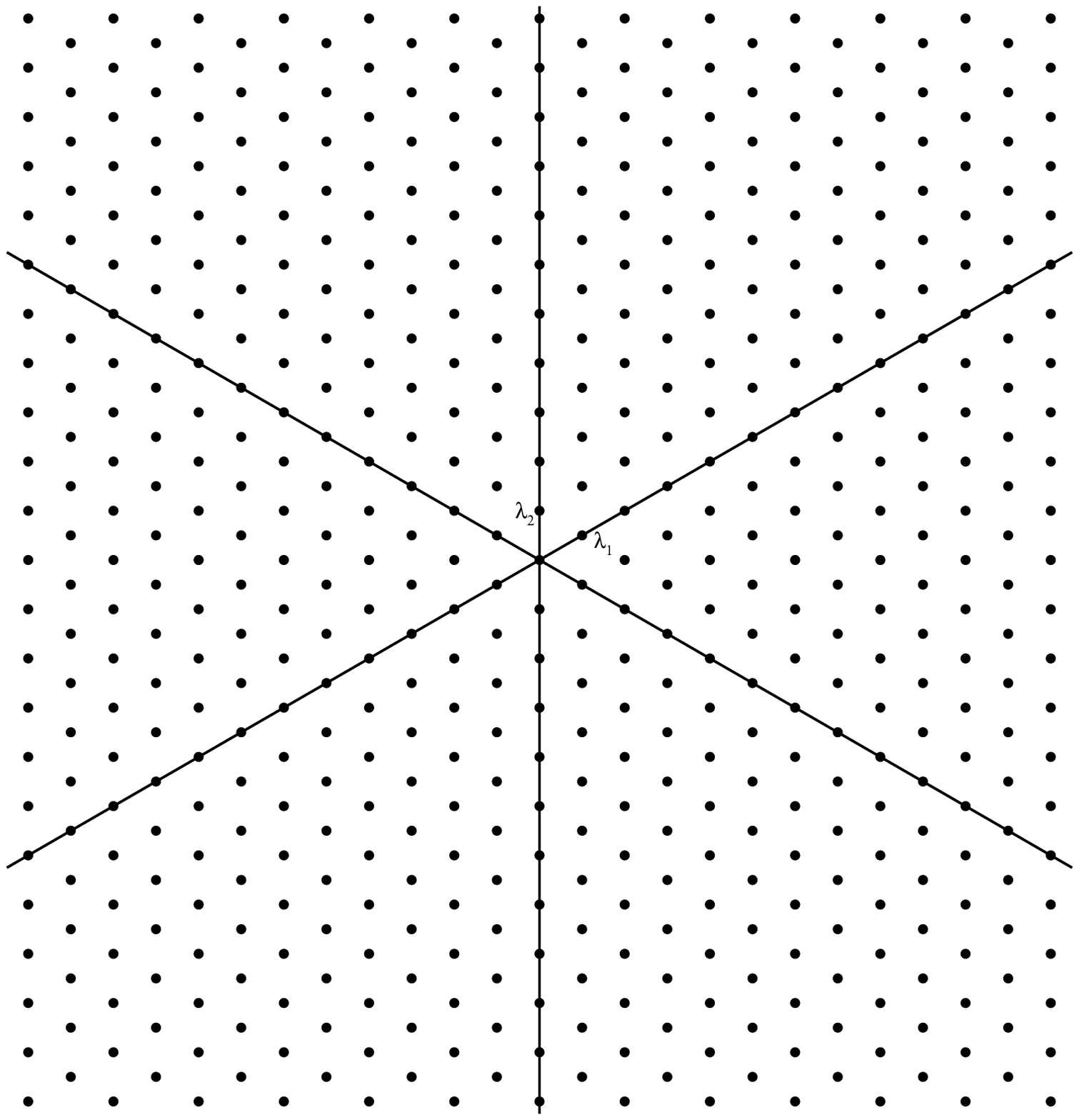,height=6.3in,width=6.5in}}

\newpage
\centerline{Figure 2: Weight Lattice of Type $B_2$}
\vskip 20pt
\centerline{\psfig{figure=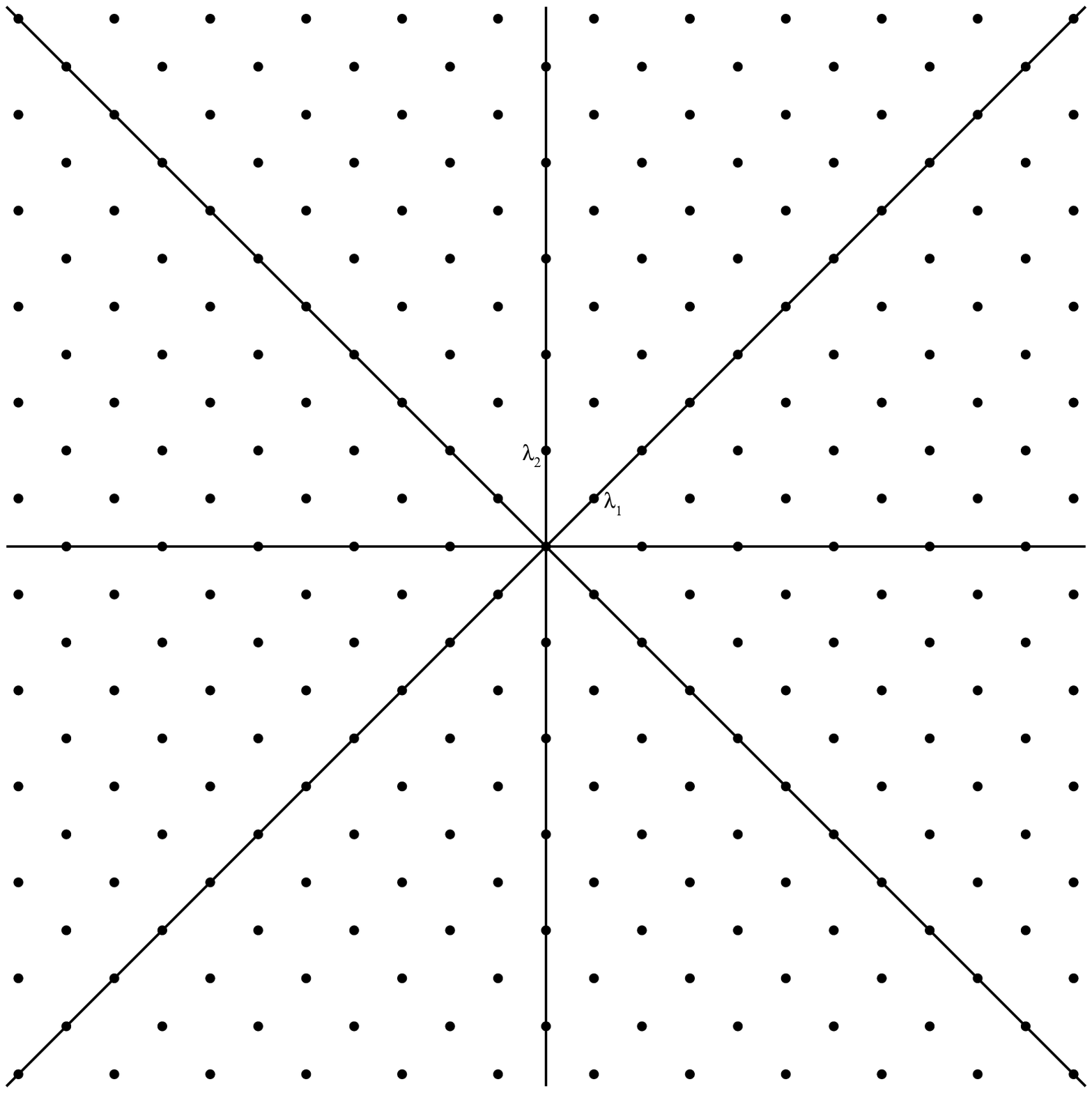,height=6.3in,width=6.5in}}

\newpage
\centerline{Figure 3: Weight Lattice of Type $G_2$}
\vskip 20pt
\centerline{\psfig{figure=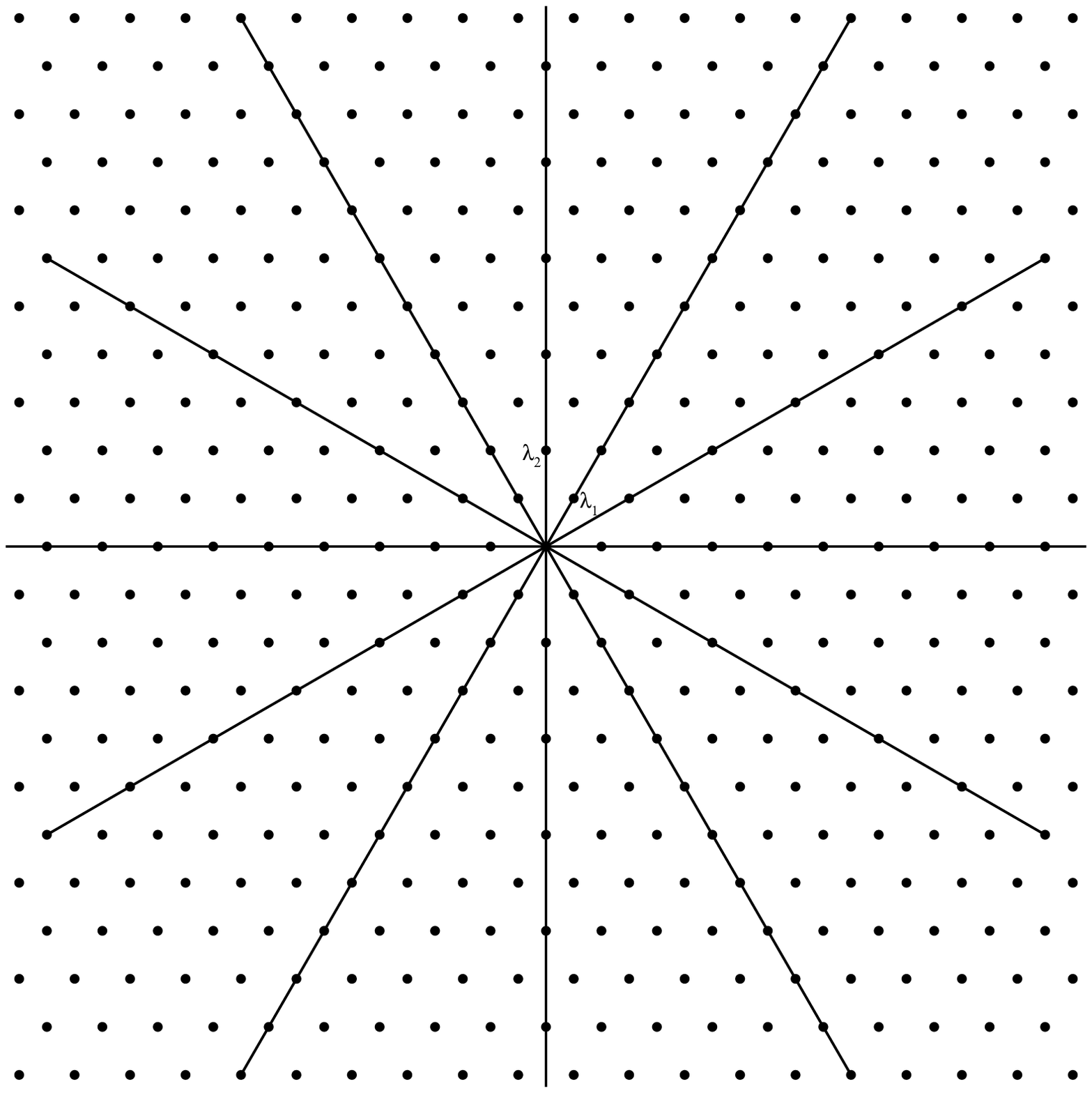,height=6.3in,width=6.5in}}

\newpage
\centerline{Figure 4: Reflection Lines of Type $A_2$}
\vskip 20pt
\centerline{\psfig{figure=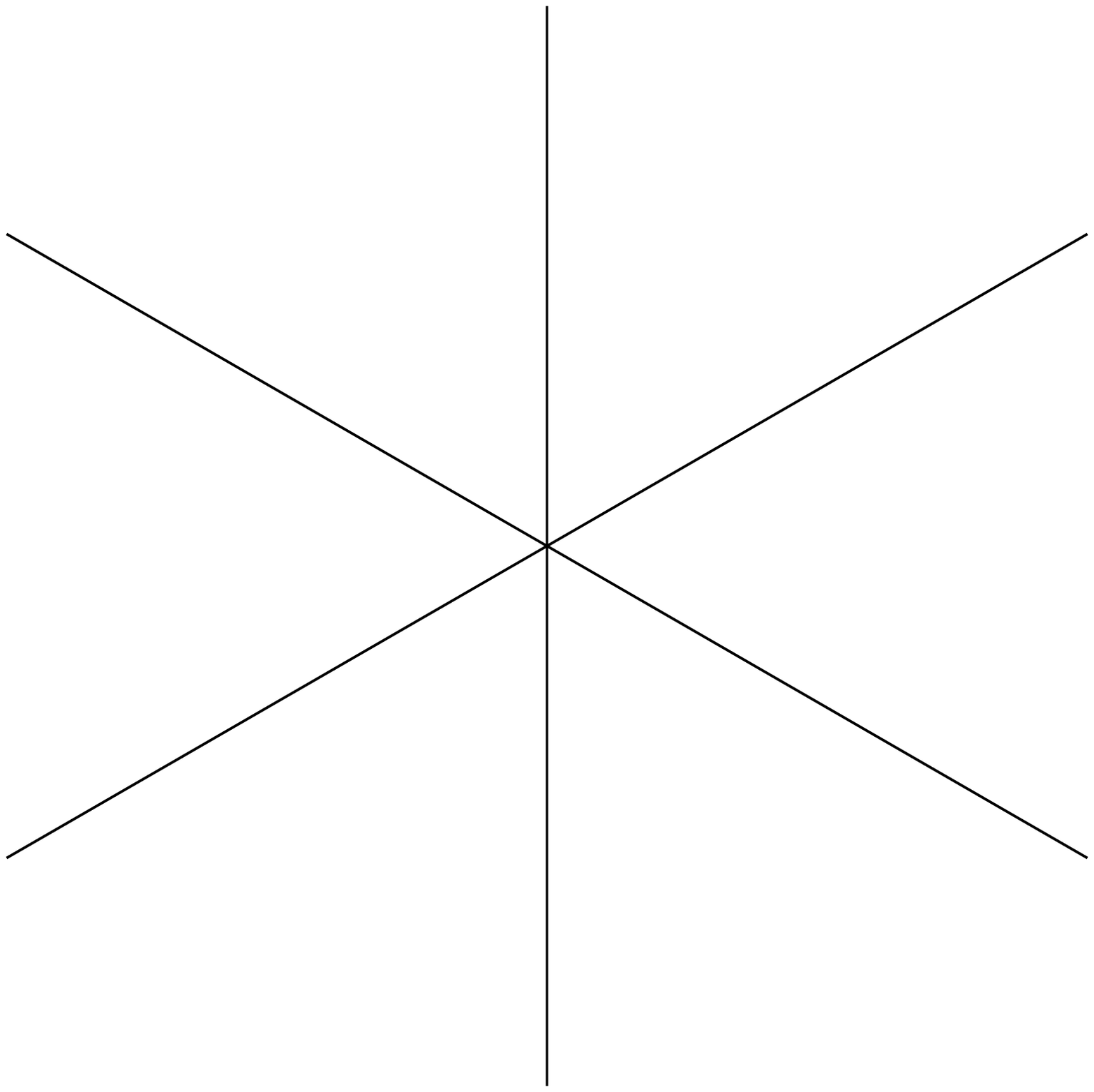,height=6.3in,width=6.5in}}

\newpage
\centerline{Figure 5: Reflection Lines of Type $B_2$}
\vskip 20pt
\centerline{\psfig{figure=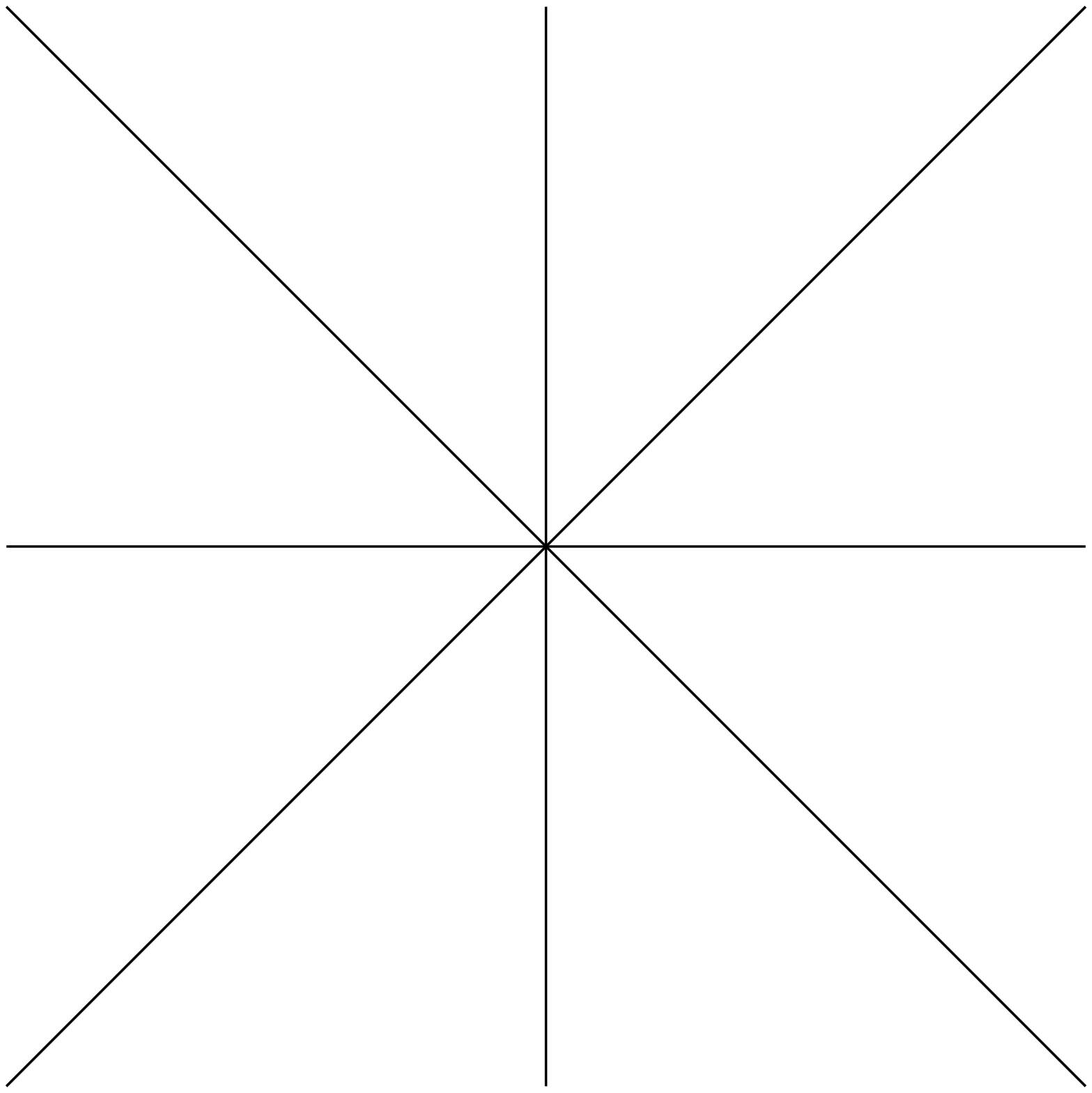,height=6.3in,width=6.5in}}

\newpage
\centerline{Figure 6: Reflection Lines of Type $G_2$}
\vskip 20pt
\centerline{\psfig{figure=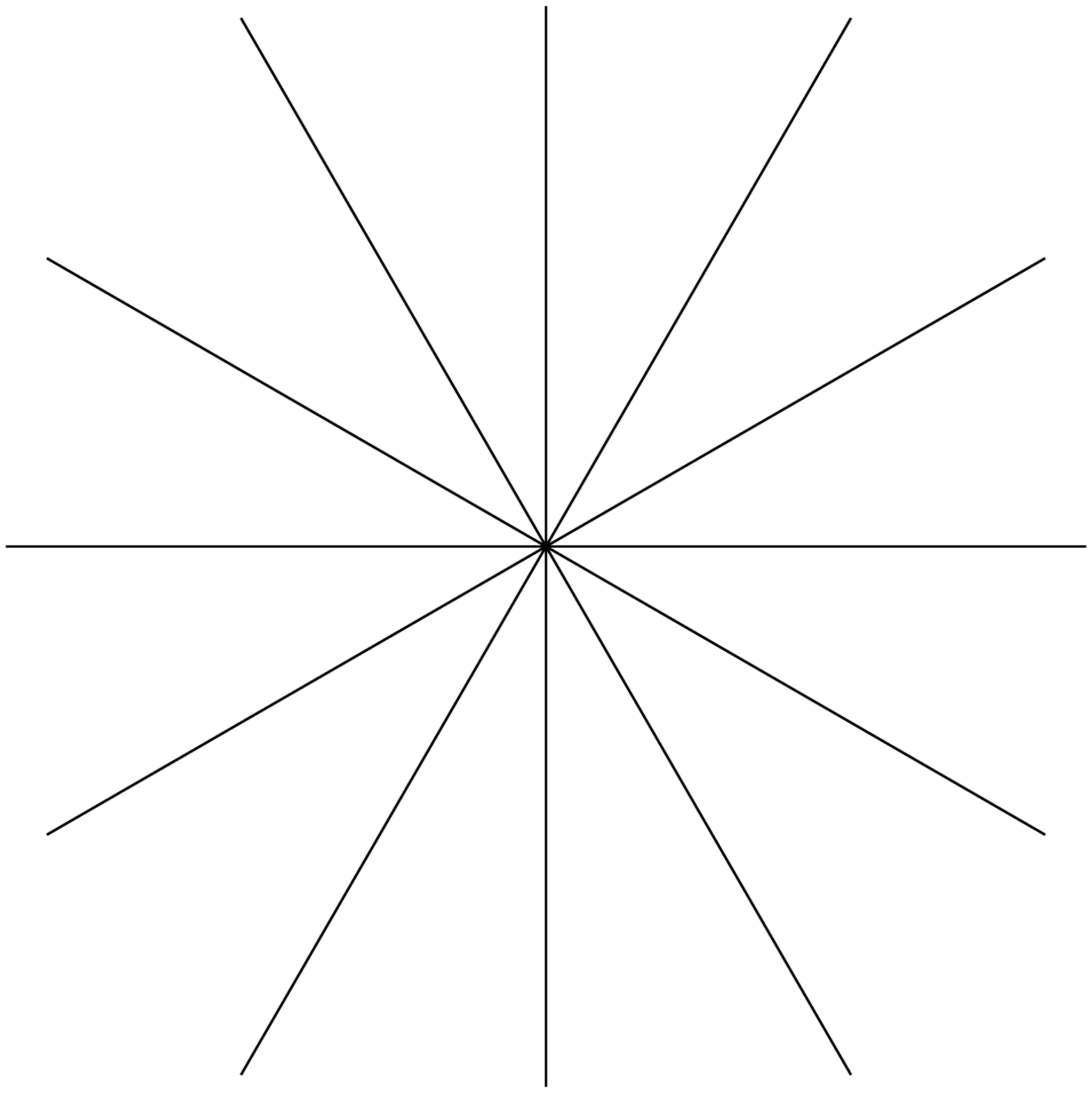,height=6.3in,width=6.5in}}

\newpage
\centerline{Figure 7: Weyl Conjugates of $\rho$ To Use In Racah Recursion}
\vskip 20pt
\centerline{\psfig{figure=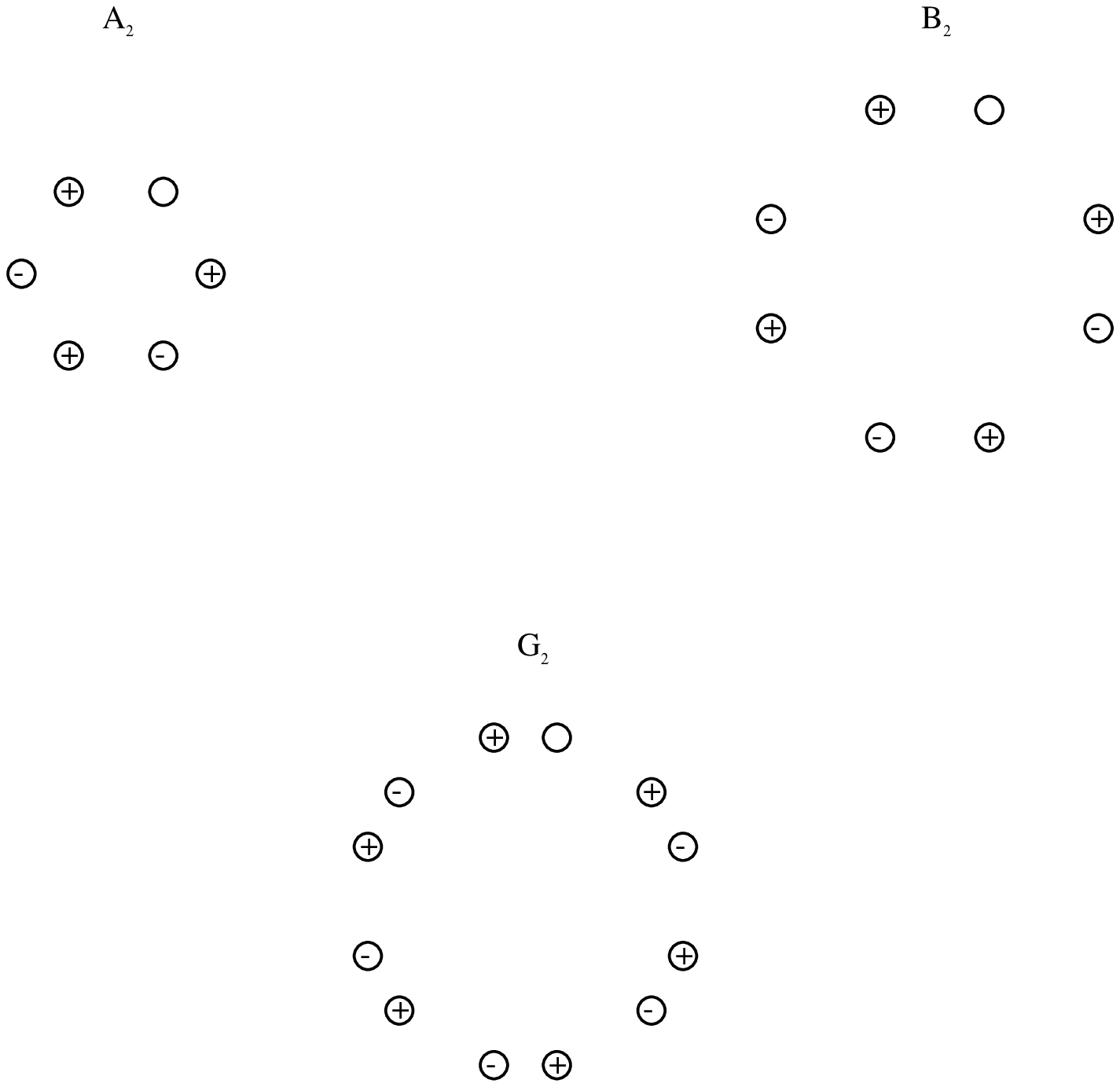,height=6.3in,width=6.5in}}

\newpage
\centerline{Figure 8: $A_2$ Weight Diagram For Irreducible Module}
\centerline{With Highest Weight $3\lambda_1 + 2\lambda_2$}
\vskip 20pt
\centerline{\psfig{figure=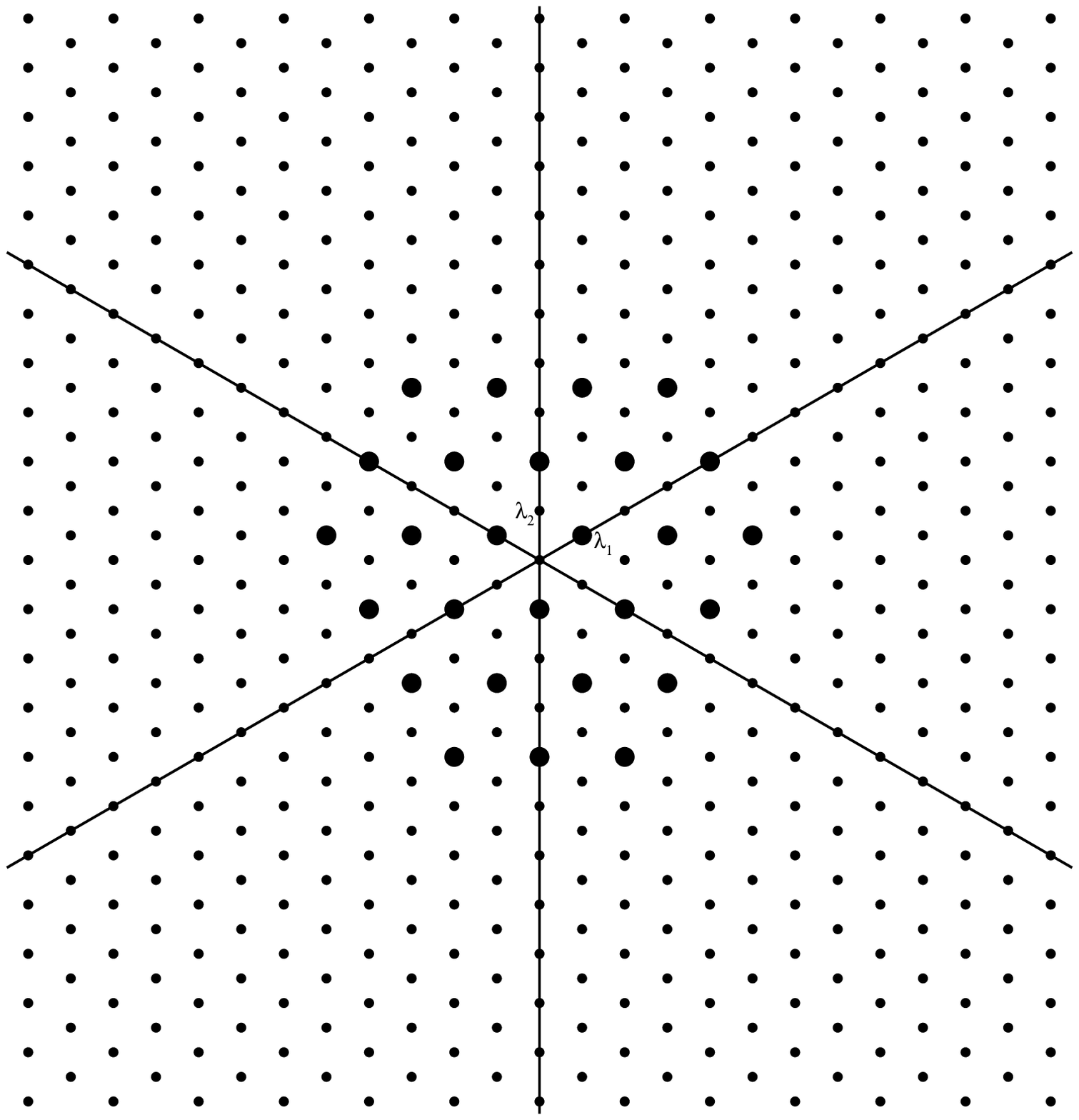,height=6.3in,width=6.5in}}

\newpage
\centerline{Figure 9: $A_2$ Weight Diagram For Irreducible Module}
\centerline{With Highest Weight $3\lambda_1 + 2\lambda_2$ Shifted by 
$\mu + \rho = 2\lambda_1 + 1\lambda_2$}
\vskip 20pt
\centerline{\psfig{figure=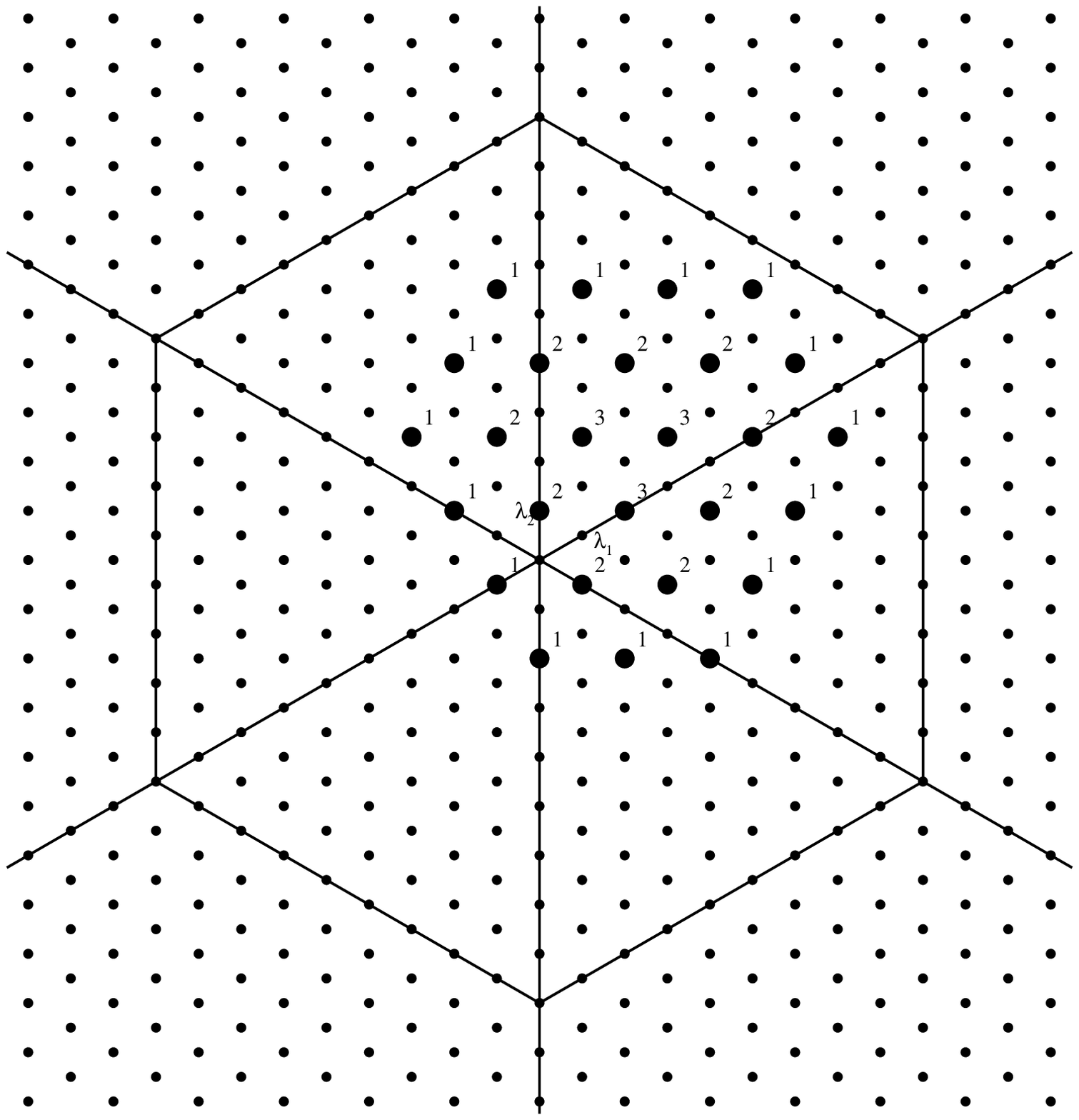,height=6.3in,width=6.5in}}

\newpage
\centerline{Figure 10: $A_2$ Weight Diagram Shifted For The}
\centerline{Level 2 Fusion Rule Computation $[3]\cdot[3] = [0]+[3]$}
\vskip 20pt
\centerline{\psfig{figure=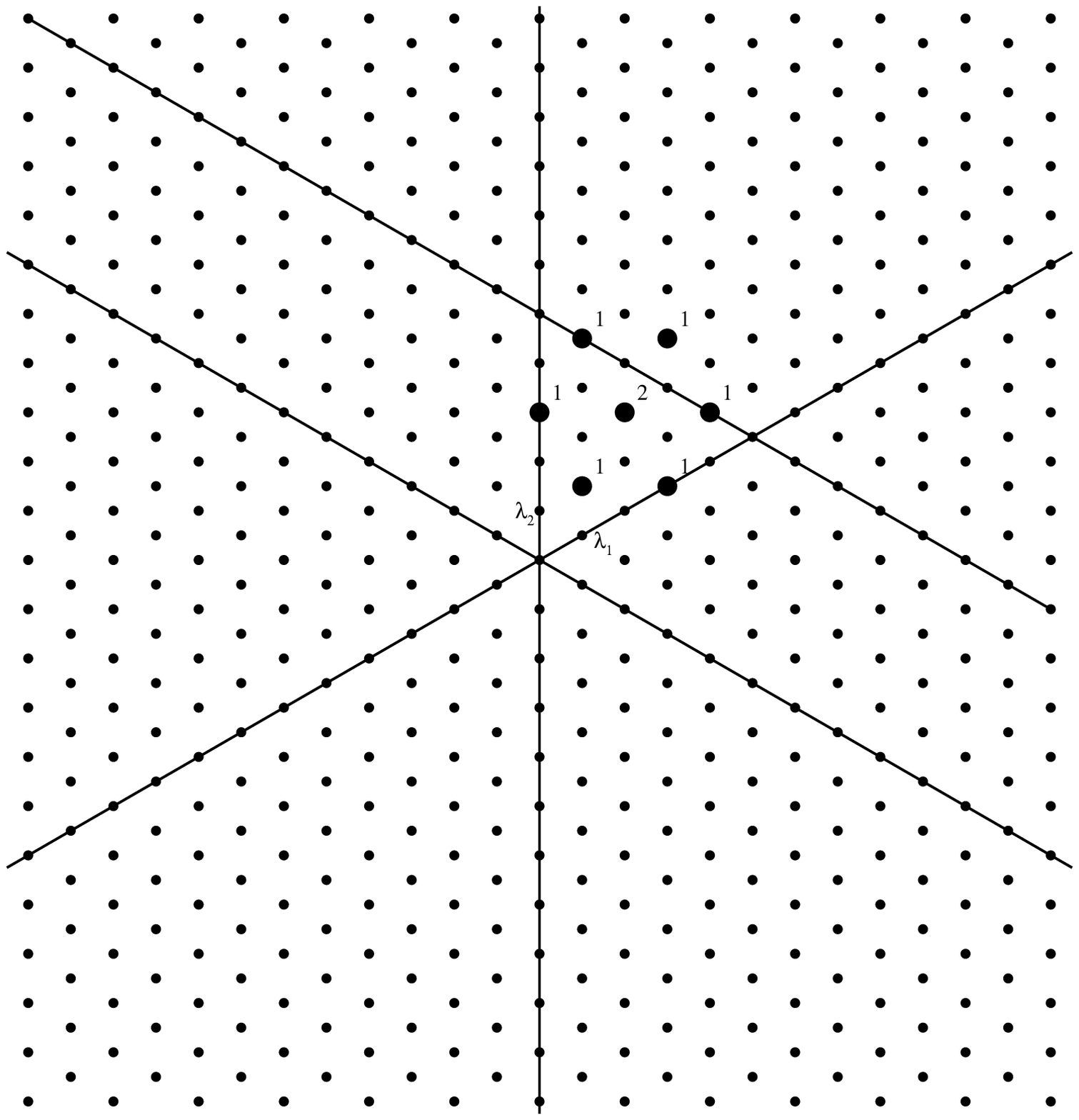,height=6.3in,width=6.5in}}

\newpage
\centerline{Figure 11: $B_2$ Weight Diagram Shifted For The}
\centerline{Level 1 Fusion Rule Computation $[1]\cdot[1] = [0]$}
\vskip 20pt
\centerline{\psfig{figure=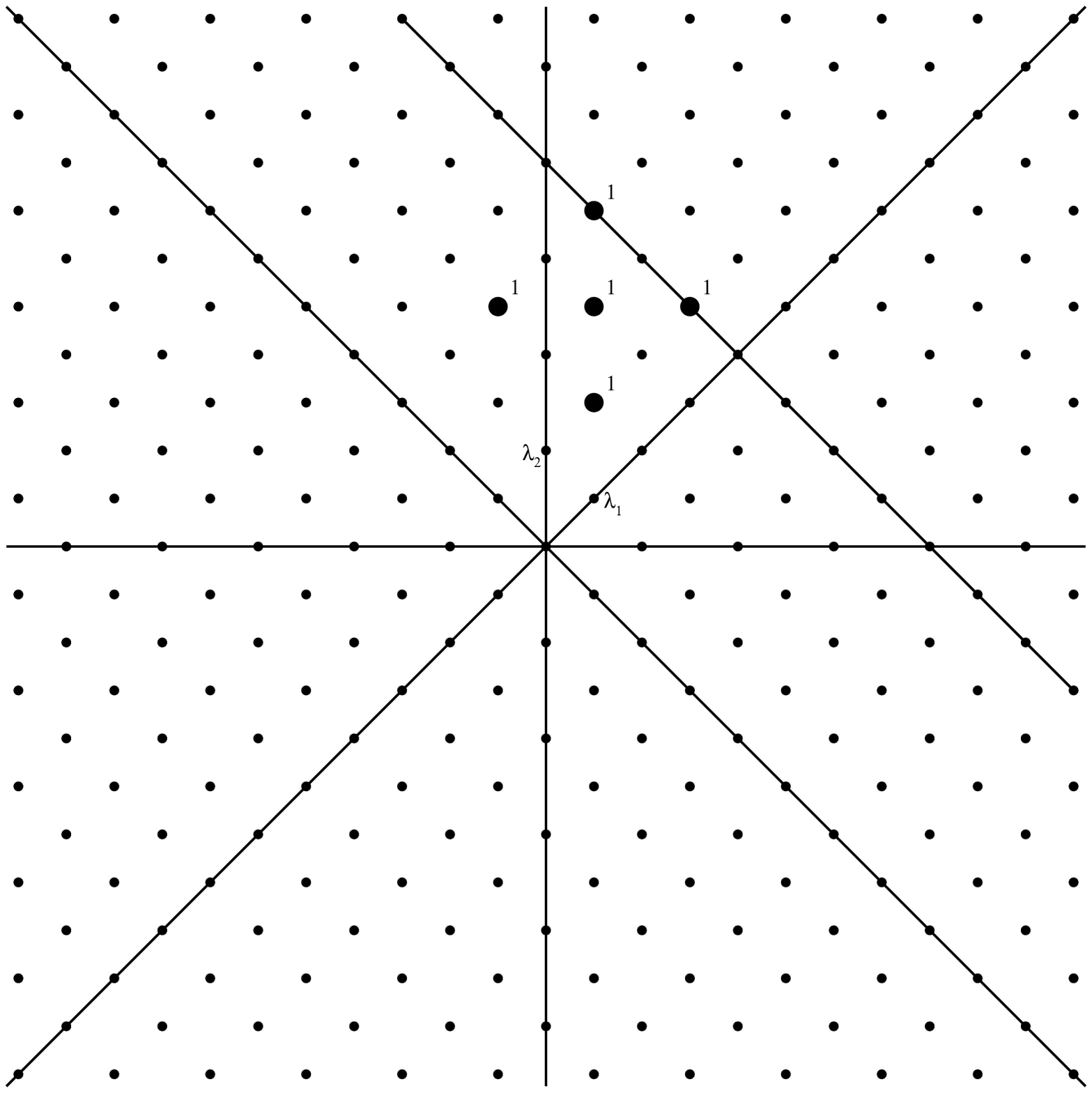,height=6.3in,width=6.5in}}

\newpage
\centerline{Figure 12: $B_2$ Weight Diagram Shifted For The}
\centerline{Level 1 Fusion Rule Computation $[1]\cdot[2] = [2]$}
\vskip 20pt
\centerline{\psfig{figure=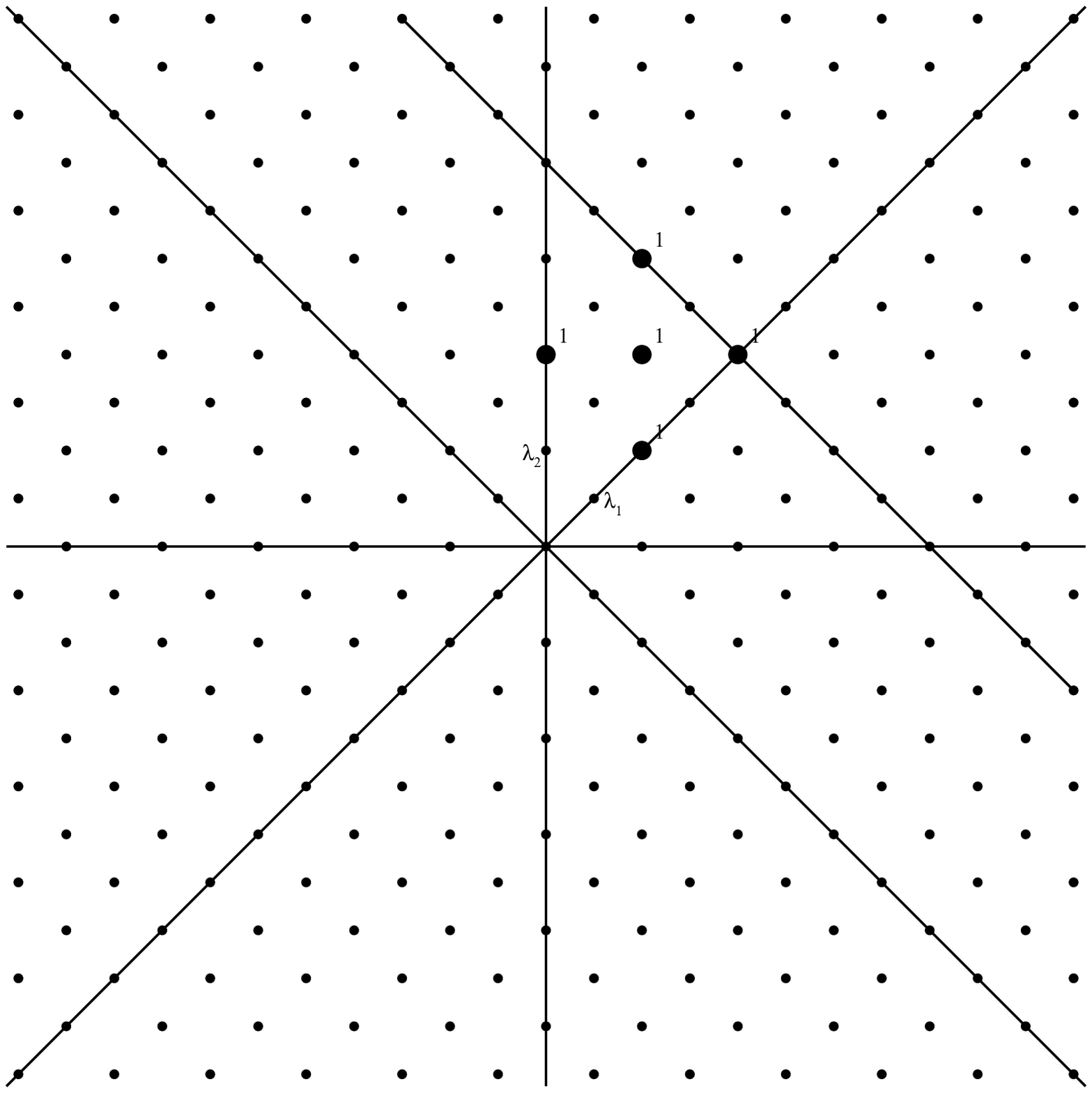,height=6.3in,width=6.5in}}

\newpage
\centerline{Figure 13: $B_2$ Weight Diagram Shifted For The}
\centerline{Level 1 Fusion Rule Computation $[2]\cdot[2] = [0]+[1]$}
\vskip 20pt
\centerline{\psfig{figure=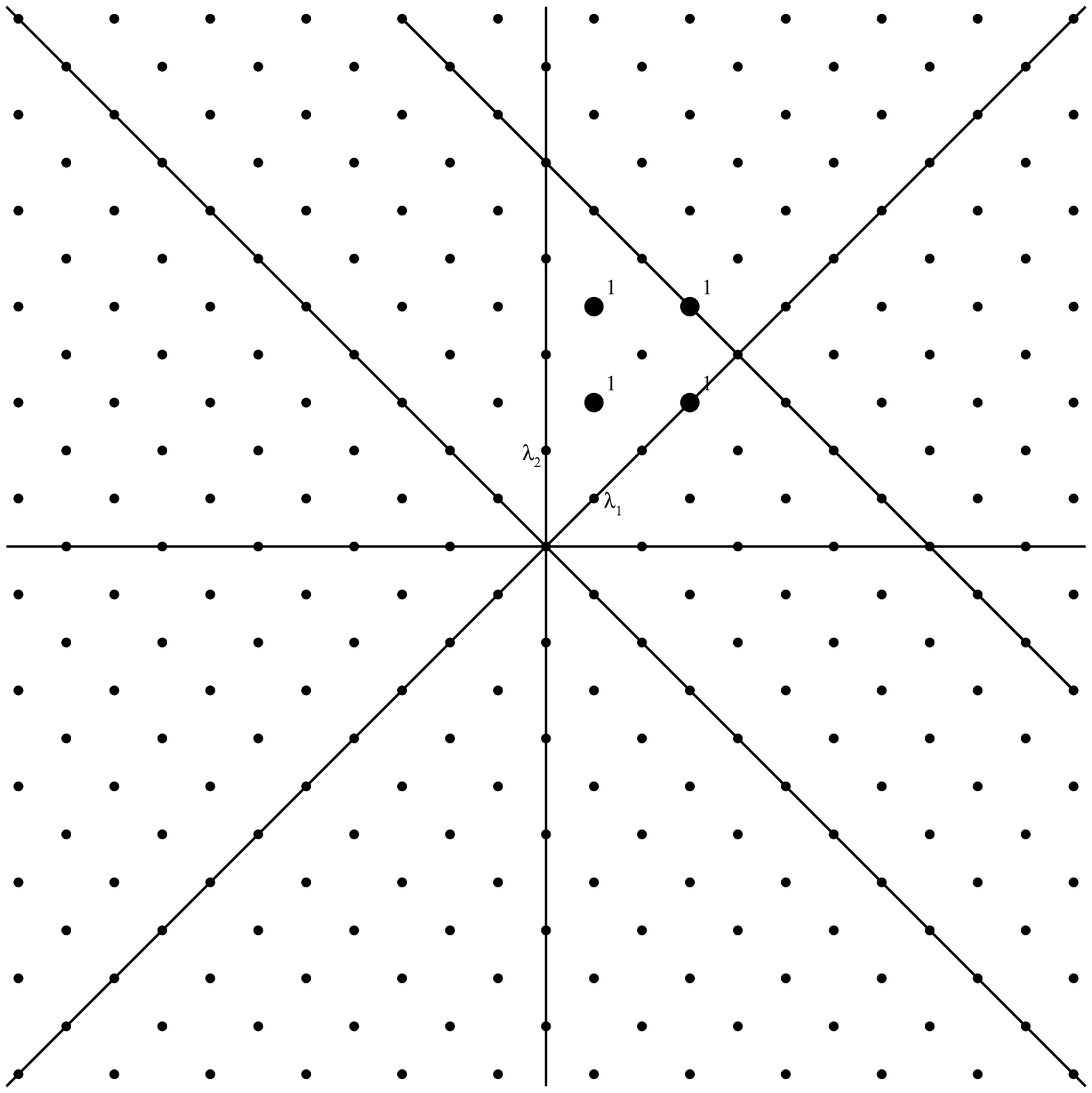,height=6.3in,width=6.5in}}

\newpage
\centerline{Figure 14: $B_2$ Weight Diagram Shifted For The}
\centerline{Level 2 Fusion Rule Computation $[3]\cdot[3] = [0]+[4]+[5]$}
\vskip 20pt
\centerline{\psfig{figure=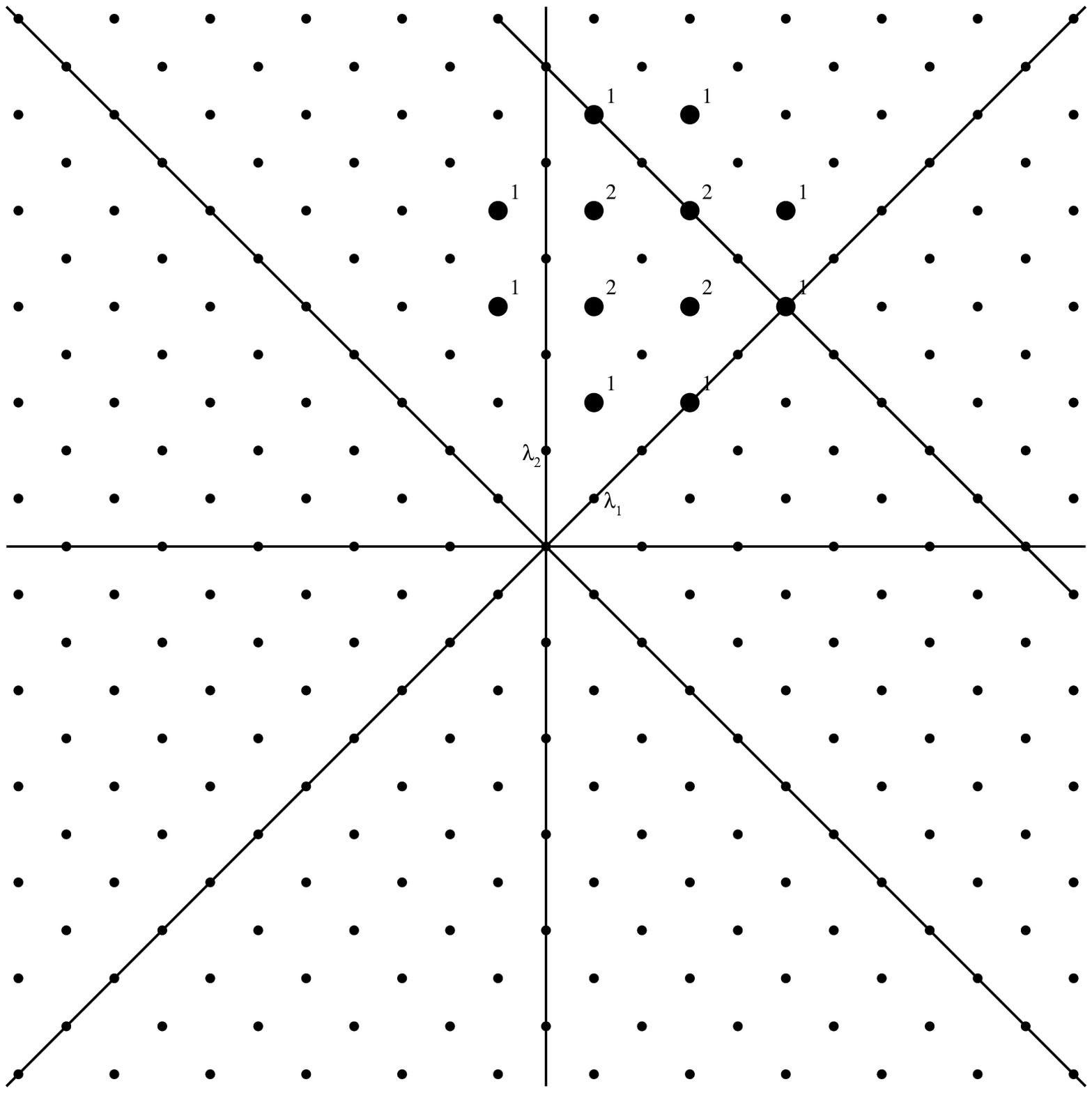,height=6.3in,width=6.5in}}

\newpage
\centerline{Figure 15: $A_2$ Weight Diagram For Irreducible Module With Highest Weight}
\centerline{$3\lambda_1 + 2\lambda_2$ Shifted by $\lambda_1+\rho$ 
For Level 5 Fusion Rule Computations}
\vskip 20pt
\centerline{\psfig{figure=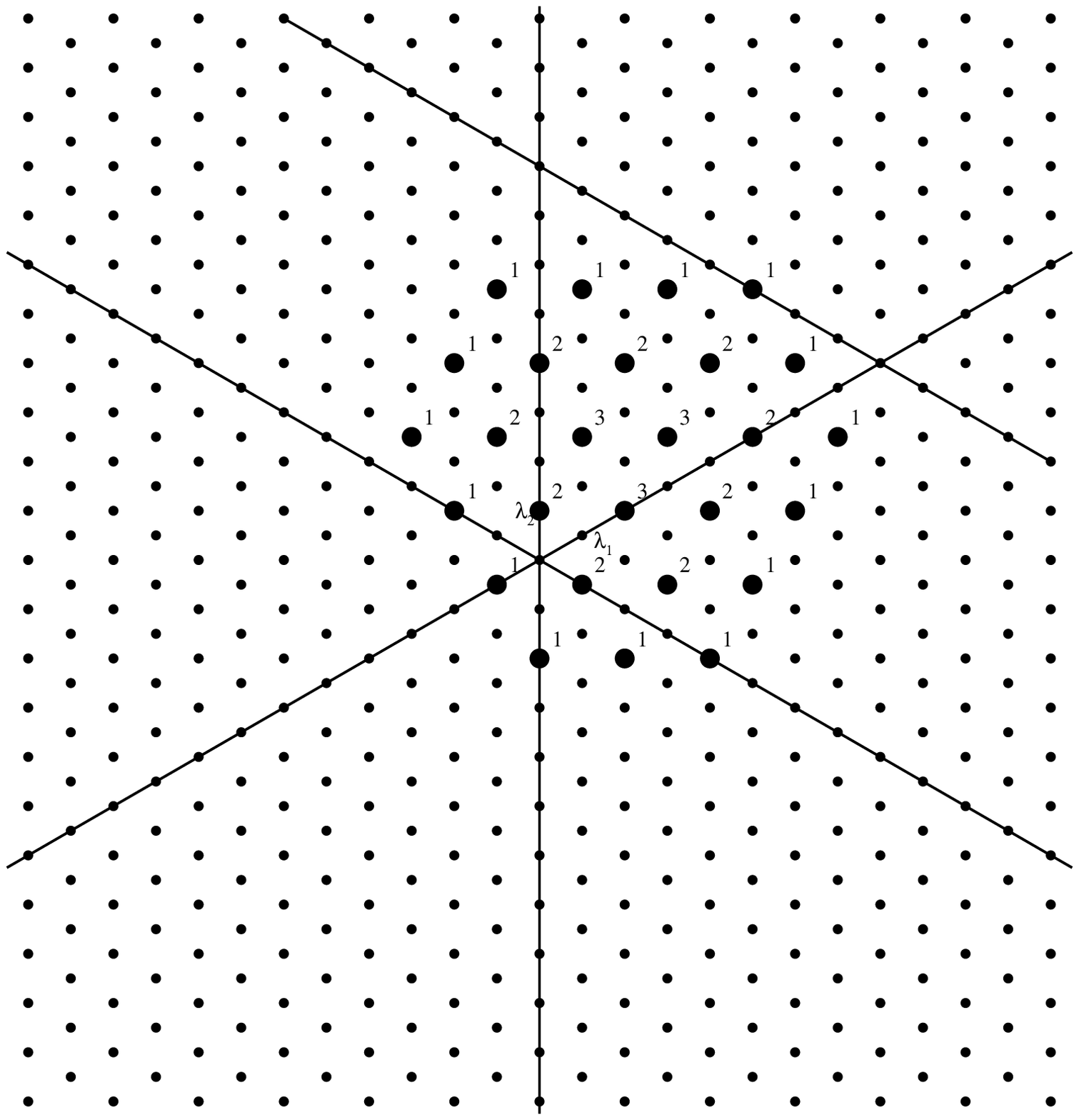,height=6.3in,width=6.5in}}

\newpage
\centerline{Figure 16: $A_2$ Weight Diagram For Irreducible Module With Highest Weight}
\centerline{$\lambda_1$ Shifted by $3\lambda_1+2\lambda_2$
 For Level 5 Fusion Rule Computations}
\vskip 20pt
\centerline{\psfig{figure=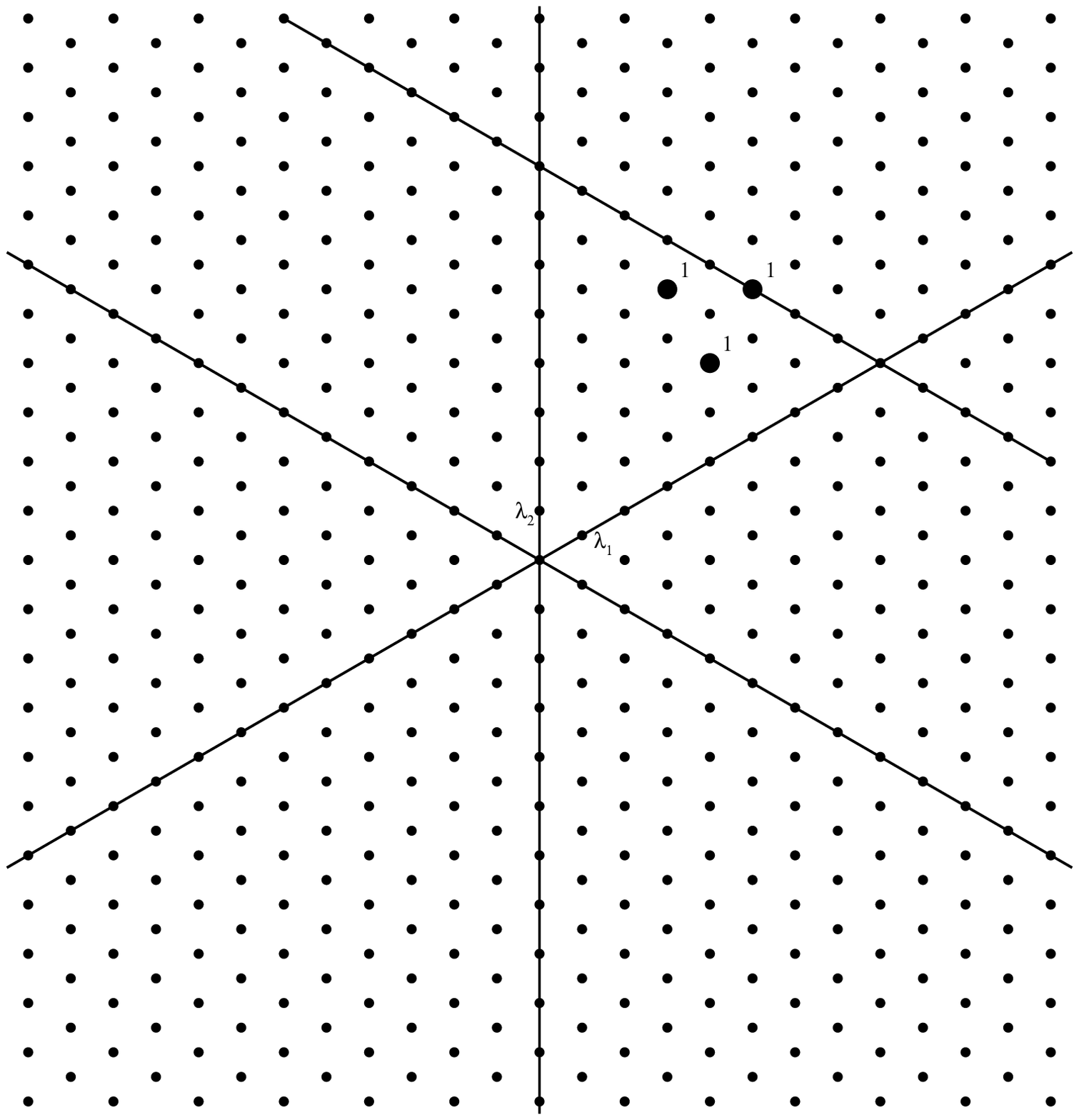,height=6.3in,width=6.5in}}

\newpage
\centerline{Figure 17: $A_2$ Weight Diagram For Irreducible Module With Highest Weight}
\centerline{$3\lambda_1 + 2\lambda_2$ Shifted by $5\lambda_1+\rho$ 
For Level 5 Fusion Rule Computations}
\vskip 20pt
\centerline{\psfig{figure=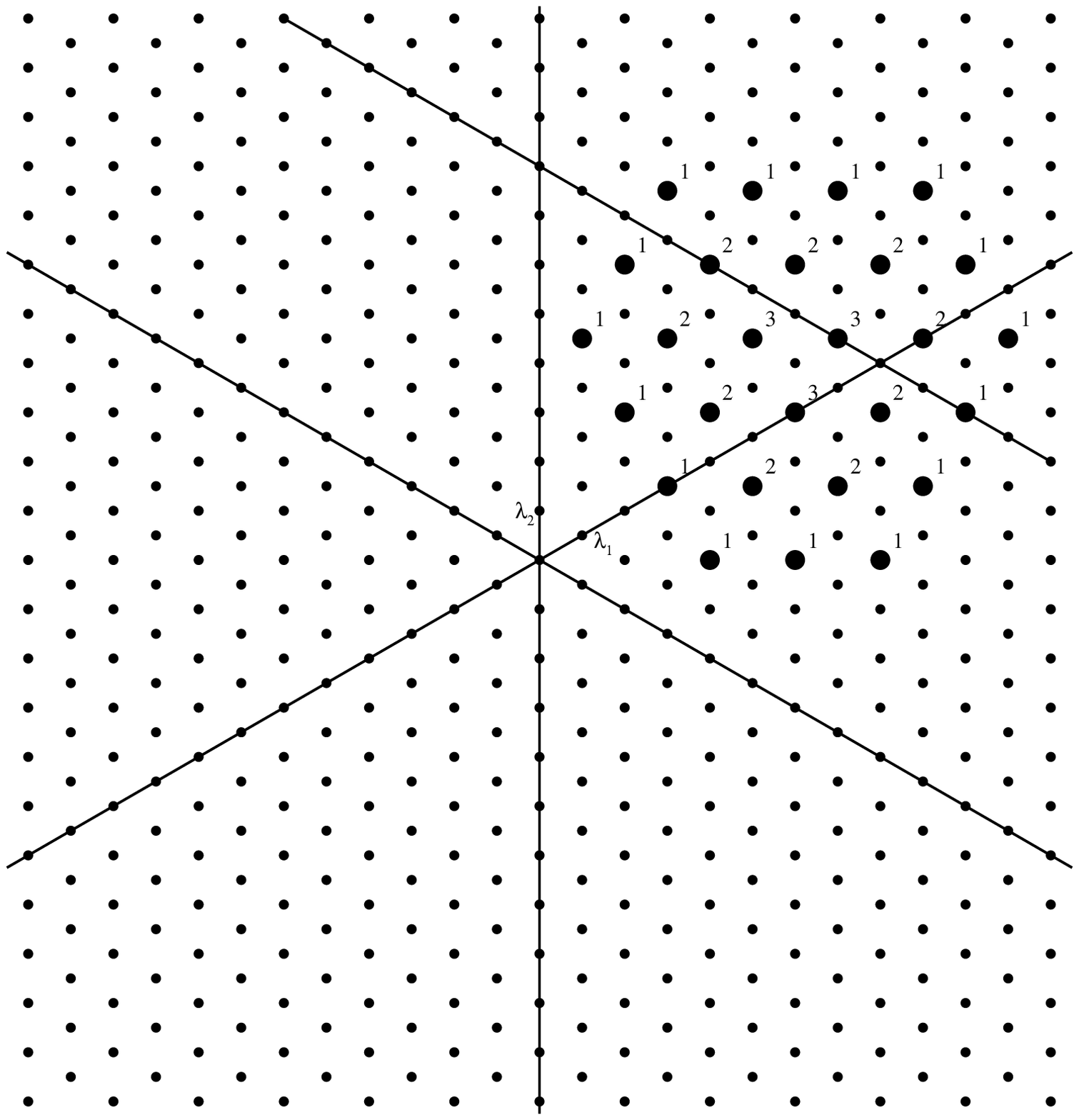,height=6.3in,width=6.5in}}

\newpage
\centerline{Figure 18: $A_2$ Weight Diagram For Irreducible Module With Highest Weight}
\centerline{$3\lambda_1 + 2\lambda_2$ Shifted by $2\lambda_1+\rho$ 
For Level 5 Fusion Rule Computations}
\vskip 20pt
\centerline{\psfig{figure=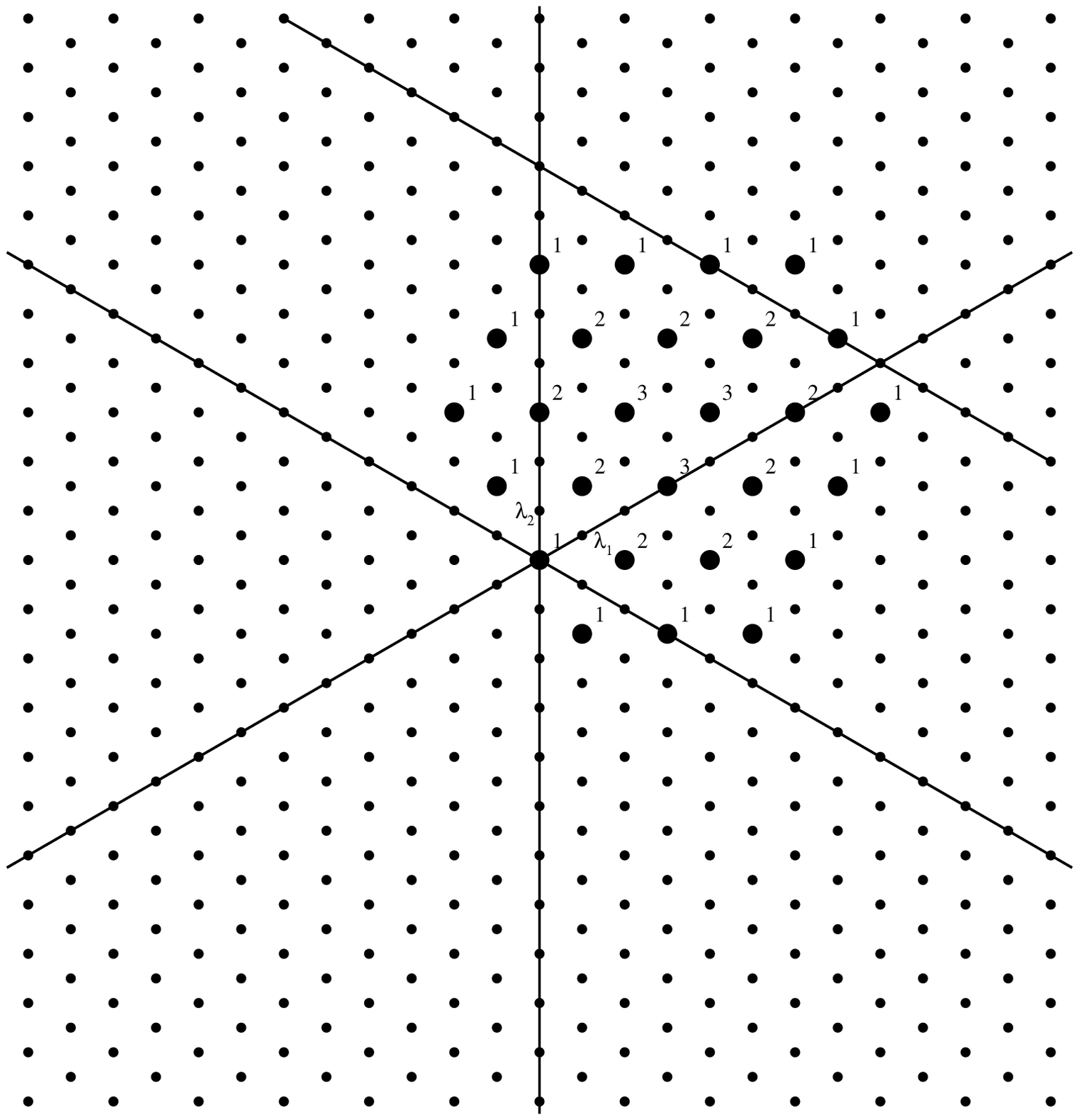,height=6.3in,width=6.5in}}

\end{document}